\documentclass[a4paper;11pt]{amsart}
\usepackage{mathrsfs}
 \usepackage[arrow,matrix]{xy}
\usepackage{amsmath,amssymb,amscd,bbm,amsthm,mathrsfs}
\newtheorem{thm}{Theorem}[section]
\newtheorem{lem}{Lemma}[section]
\newtheorem{cor}{Corollary}[section]
\newtheorem{prop}{Proposition}[section]
\newtheorem{rem}{Remark}[section]

\theoremstyle{definition}

\begin{document}
\numberwithin{equation}{section}

\title[ On   pluripotential theory associated to quaternionic $m$-subharmonic functions]{ On     pluripotential theory associated to quaternionic $m$-subharmonic functions}
\author{Shengqiu Liu and Wei Wang}
\begin{abstract}Many aspects of pluripotential theory are generalized to quaternionic $m$-subharmonic functions. We introduce quaternionic version of notions of the $m$-Hessian operator,   $m$-subharmonic functions, $m$-Hessian measure,  $m$-capapcity, the relative $m$-extremal function and
the $m$-Lelong number,  and show various propositions for them, based on $d_0$ and $ d_1$ operators, the  quaternionic counterpart of $\partial$ and $\overline{\partial}$, and quaternionic closed positve currents.   The definition of quaternionic $m$-Hessian operator can be extended to locally bounded quaternionic $m$-subharmonic functions and   the corresponding convergence theorem is proved.
  The  comparison principle and the quasicontinuity of   quaternionic $m$-subharmonic functions are established.
We also find the fundamental solution  of the quaternionic  $m$-Hessian operator.
\end{abstract}\thanks{
Supported by National Nature Science Foundation in China (No.
11971425)  }\thanks{   Department of Mathematics,
Zhejiang University, Zhejiang 310027,
 P. R. China, Email:  liushengqiu97@163.com, wwang@zju.edu.cn}
\maketitle
\section{\textbf{Introduction}}

Pluripotential theory    provides fine properties of plurisubharmonic   functions, their Monge-Amp\`ere measure and solutions to the complex  Monge-Amp\`ere equation $(dd^c u)^n=f\beta^{n }$, where $\beta $ is the fundamental K\"ahler form on $\mathbb{C}^n$. Notably the Monge-Amp\`ere operator $(dd^c u)^n $ is well defined for some non-smooth plurisubharmonic functions,
e.g. continuous or locally bounded  plurisubharmonic functions. This theory  is a powerful tool in complex analysis  of several variables, and was generalized to $m$-subharmonic functions, their Hessian measure and the complex $m$-Hessian equation $(dd^c u)^m\wedge\beta^{n-m}=f\beta^{n }$. Pluripotential theory for $m$-subharmonic functions developed rapidly in last two decades, and there are vast literatures (cf. \cite{ACH,AC,Benali,blocki,dinew3,Elkhadhra,Hai,Hung,li,Lu,Nguyen,NguyenK,sadullaev1,sadullaev2,wan-wang1} and references therein).

 On the quaternionic space,   Alesker  \cite{alesker1} introduced   notions of   quaternionic plurisubharmonic functions and quaternionic Monge-Amp\`ere  operator,  proved a quaternionic version of the Chern-Levine-Nirenberg estimate
and extended  the quaternionic Monge-Amp\`ere  operator to continuous
quaternionic plurisubharmonic functions.
He also \cite{alesker2} used the Baston operator $\triangle$ to express the quaternionic Monge-Amp\`{e}re operator by using   methods of complex geometry. Then
Wan-Wang  \cite{wan-wang1} introduced the first-order differential operators $d_0$ and $d_1$ acting on the quaternionic version of differential forms   and the notion of the closedness of  a quaternionic positve current, motivated by  $0$-Cauchy-Fueter complex in quaternionic analysis \cite{Wang}. The behavior of $d_0,d_1$ and $\Delta=d_0d_1$ is very similar to $\partial,\overline{\partial}$ and $\partial \overline{\partial}$ in several complex variables, and many results in the complex pluripotential theory   have been also extended to the quaternionic case (cf.
\cite{alesker4,alesker6,Bou,wan19,wan20,wan5,wan-wang,WZ,wang2,wang21} and references therein). Some aspects of quaternionic  pluripotential theory has been generalized to the Heisenberg group \cite{wang21}.
  The purpose of this paper is to generalize     pluripotential theory to   quaternionic $m$-subharmonic functions.

  The paper is organized as follows. In Section 2,   a quaternionic version of Garding inequality is given by applying   Garding's theory of hyperbolic polynomials  to symmetric function  of eigenvalues of a quaternionic hyperhermitian matrix. In Section 3, we briefly recall   positive forms,  the first-order differential operators $d_0$ and $d_1$ and $\Delta=d_0d_1$ and their various  propositions.   The   quaternionic   $m$-Hessian operator   is introduced and   can   be written as $(\Delta u)^{m}\wedge  \beta_n ^{n-m}$, where $\beta_n $ is the fundamental   form on $\mathbb{H}^n$.
  In Section 4, we give the definition  of nonsmooth quaternionic $m$-subharmonic function in terms of  positive currents, which coincides with that for smooth ones, and prove basic properties of   quaternionic $m$-subharmonic functions.
    In Section 5, for continuous quaternionic $m$-subharmonic functions,  the locally uniform estimate,  i.e. the Chern-Levine-Nirenberg estimate, the existence of $m$-Hessian measure and the  comparison principle are established. We study   the relative $m$-extremal function  and quaternionic $m$-capapcity in Section 6, and establish the quasicontinuity of  quaternionic $m$-subharmonic functions, the extension of   quaternionic $m$-Hessian operator to locally bounded  quaternionic $m$-subharmonic functions and the corresponding convergence theorem (the Bedford-Taylor theory) in Section 7.
    In Section 8 we find the fundamental solution  of the   $m$-Hessian operator   and define the $m$-Lelong number for a quaternionic $m$-subharmonic function.

We use the Sadullaev-Abdullaev approach \cite{sadullaev1,sadullaev2} to  $m$-subharmonic functions and the complex $m$-Hessian operator, i.e.
based on an integral estimate for $\int_{\Omega}(\Delta u)^m\wedge\beta_n^{n-m}$ on a domain $\Omega$. While in the classical approach (e.g. \cite{klimek}),
ones usually only use local estimate by using a cut-off function,  e. g. in the proof of the Chern-Levine-Nirenberg estimate. We established such integral estimate by using a Stokes-type formula instead of Stokes  formula, since our forms are not differential forms. The advantage of
this approach is that we can quite quickly  to establish   necessary estimates and various results.

\section{\textbf{Hyperbolicity of symmetric functions of eigenvalues of a quaternionic hyperhermitian matrix}}
\subsection{Quaternionic hyperhermitian matrix}
An $n\times n$ quaternionic   matrix $A=(a_{ij})$ is called  {\it  hyperhermitian} if $A^*=A$, i.e., $a_{ij}=\overline{a}_{ji}$ for all $i,j$.  Denote by $\mathscr{H}^n $   the   space of all quaternionic hyperhermitian $n\times n$ matrices,  by $GL_{\mathbb{H}}(n)$ the set of all invertible quaternionic $(n\times n)$-matrices, and   by $U_{\mathbb{H}}(n)$ the set of all unitary quaternionic $(n\times n)$-matrices, i.e.
       $
U_{\mathbb{H}}(n)=\{{{M}}\in GL_{\mathbb{H}}(n),{{M}}^*{{M}}={{M}}{{M}}^*=I_n\}.
 $
Let us recall the definition of the Moore determinant \cite{Aslaksen}
  for ${M}=({M}_{ij})\in\mathscr{H}^n$. Write a permutation $\sigma$ of $(1,\dots,n)$ as a product of disjoint cycles as
\begin{equation*}
\begin{aligned}
\sigma=(n_{11}\dots n_{1l_1})(n_{21} \cdots n_{2l_2})\cdots(n_{r1}\cdots n_{rl_r}),
\end{aligned}
\end{equation*}
where for each $i$, we have $n_{i1}<n_{ij}$ for all $j>1$, and $n_{11}>\dots>n_{r1}$. Then   \begin{equation}
\begin{aligned}\label{moore determiniant}
\det M=\sum\limits_{\sigma\in S_n}\text{sgn} \sigma {M}_{n_{11}n_{12}}\cdots {M}_{n_{1l_1}n_{11}}M_{n_{21}M_{22}}\cdots {M}_{n_{rl_r}n_{r1}}.
\end{aligned}
\end{equation}
 Consider the homogeneous polynomial $\text{det}(s_1{{M}}_1+\ldots+s_n{{M}}_n)$ in real variables $s_1,\ldots,s_n$ of degree $n$. The coefficient of the monomial $s_1\cdots s_n$ divided by $n!$ is called the {\it mixed determiniant} of the hyperhermitian matrices ${{M}}_1,\ldots,{{M}}_n$, and   is denoted by $\text{det}({{M}}_1,\ldots,{{M}}_n)$.

  \begin{prop} \label{diagonal and real}
 (1) \cite[Claim 1.1.4, 1.1.7]{alesker1}  For a hyperhermitian $(n\times n)$-matrix ${{M}}$, there exits a unitary ${U}$ such that ${U}^*{{M}U}$ is diagonal and real. \\
  (2) \cite[Theorem 1.1.9]{alesker1}
  for any quaternionic hyperhermitian $(n\times n)$-matrix $ {M}$ and for any quaternionic $(n\times n)$-matrix $ {C}$, we have
$
\det({{ C}}^*{{M}}{C})=\det({{M}})\det({{ C}}^*{{ C}}).
$
\\
 (3)  \cite[P. 11]{alesker1}
  The mixed determinant is symmetric with respect to all variables, and linear with respect to each of them. In particular, $\det({A},\dots,{A})=\det({A})$.
\end{prop}

   \subsection{Hyperbolic polynomials} Recall  Garding's  theory of hyperbolic polynomials  \cite{garding}. Let $P $ be a homogeneous polynomial of degree $m $ in  variables $x \in\mathbb{R}^N$. We say that $P$ is {\it hyperbolic at  $a\in\mathbb{R}^N$}  if the equation $P(sa + x) = 0$ has $m$ real zeros for every $x\in\mathbb{R}^N$.
The {\it completely polarized form} of the polynomial $P$ is given by
\begin{equation}\label{polar}
M(x^1,\dots,x^m)=\frac{1}{m!}\prod_k\left(\sum\limits_i x_i^k\frac{\partial}{\partial x_i}\right)P(x),\end{equation}
where $x^k=(x_1^k,\dots,x_N^k),x=(x_1,\dots,x_N)\in\mathbb{R}^N.$

Let $ C(P,a)$ be the set of all $x\in \mathbb{R}^N$ such that   $P(sa+x)\not= 0$ when $s\geq 0$.
If we factorize it as
$
P(sa+x)=P(a)\prod_{1}^m(s+\mu_k(a,x)),
$
for fixed   $x\in\mathbb{R}^N $,  then $x\in C(P,a)$ is equivalent to require
\begin{equation}\label{eq:h(a,x)}
h(a,x): =\min_k\mu_k(a,x) >0.
\end{equation}
  The {\it{linearlity}} $LP$ of $P$ is defined as the set of all $x$ such that $P(sx+y)=P(y)$ for all $s$ and $y$. The edge $\partial C$ of $C=C(P,a)$ is the set of all $x$ such that $C+x=C$ (cf.  \cite[ P. 962]{garding}).
\begin{prop}
Suppose a homogeneous polynomial $P$ on $\mathbb{R}^N$ of degree $m>1$   is hyperbolic at $a\in\mathbb{R}^N$. Then
  (1)\label{prop:derivatives} \cite[Lemma 1]{garding}
$
Q =\sum\limits_{k=1}^N a_k\frac{\partial P}{\partial x_k}
$
is hyperbolic at $ a$.

 (2) \cite[Theorem 2]{garding}
The function $h$ defined in  \eqref{eq:h(a,x)}  is positive, homogeneous and concave, i.e.
$
h(a,sx)=sh(a,x)$ for $s\geq 0 $ and $ h(a,x+y)\geq h(a,x)+h(a,y)$.
In particular, $C=C(P,a)$ is convex. Further, $P$ is hyperbolic at any $b \in C$ and $C(P,b)=C(P,a).$

 (3)  \cite[Theorem 3]{garding}
  $\partial C=LP$ and $x$ belongs to $LP$ if and only if
$
\mu_1(a,x)=\dots=\mu_m(a,x)=0.
$
\end{prop}
\begin{prop} \label{point}
{\rm \cite[Theorem 5]{garding} }
Let a homogeneous polynomial $P$ of degree $m>1$ be hyperbolic at $a\in\mathbb{R}^N$, $P(a)>0$   and let $M$ be the completely polarized form of $P$. If $x^1,\dots,x^m\in  C(P,a)$, then
\begin{equation}
M(x^1,\dots,x^m)\geq P(x^1)^{\frac{1}{m}}\dots P(x^m)^{\frac{1}{m}}
\end{equation}
with equality if and only if $x^1,\dots,x^m$ are pairwise proportional modulo LP.
\end{prop}
\subsection{The hyperbolicity of symmetric functions of eigenvalues of a quaternionic hyperhermitian matrix}
\rm Now we apply the above theory of hyperbolic polynomials to symmetric functions of eigenvalues of a quaternionic hyperhermitian matrix.
An element $x=(x_{ij})\in\mathscr{H}^n$ is
$1$-$1$ correspondent to a point   $(x_{12},\dots,x_{(n-1)n},x_{11},\dots,x_{nn})$ in $\mathscr{H}^{\frac{n(n-1)}{2}}\times \mathbb{R}^n$. So we can identify $\mathscr{H}^n$ with $\mathbb{R}^N$ for $N=2n^2-n$.
  \begin{prop}\label{prop: hyperbolic'}
  $
P(x)=\det x$
is hyperbolic   at $I$   on $\mathscr{H}^n$, where $I$ is the identity matrix in $\mathscr{H}^n$.
\end{prop}
\begin{proof}
By definition $(\ref{moore determiniant})$ of the Moore determinant,   we can write $\det x=Q_1(x)+\mathbf{i}Q_2(x)+\mathbf{j}Q_3(x)+\mathbf{k}Q_4(x)$ for some real polynomials $Q_1,\dots,Q_4$ of degree $n$. On the other hand by Proposition  \ref{diagonal and real} (1), we have  $\det x=\prod\limits_{k}\lambda_k(x)\in\mathbb{R}$ with  $\lambda_k(x)$ ($ k=1,\dots,n)$  to be   eigenvalues of the hyperhermitian matrix $x$, which are all real. We see that
$
\det x=Q_1(x)
$.
So $P(x)=\det x$ is a real polynomial of degree $n  $.
It follows from Proposition $\ref{diagonal and real}$ (2) that there exists a unitary matrix $U$ such that $x={U}\operatorname {diag}(\lambda_1,\dots,\lambda_n){U}^*,$ and so
\begin{equation}
\begin{aligned}\label {AAAAA}
P(sI+x)=\det(sI+\operatorname {diag} (\lambda_1 ,\dots,\lambda_n) )=\prod\limits_1^n(s+\lambda_k ).
\end{aligned}
\end{equation}
Therefore $P(sI+x)$ has exactly $n$ real zeros, i.e. $P(x)=\det x$ is hyperbolic at $I$.
\end{proof}

For $A \in\mathscr{H}^n$, let $\lambda_1(A)\leq  \dots\leq  \lambda_n(A)$ be eigenvalues of $A$ and write $\lambda(A)=(\lambda_1(A),\dots,\lambda_n(A))$   as a vector in $\mathbb{R}^n$. Set
      \begin{equation}\label{EEEEE}
 \mathcal{ H}_m(A):=S_m(\lambda(A)),\end{equation}
 where
\begin{equation}\label{eq:Sm}
   S_m(\lambda)=\sum\limits_{1\leq j_1<\cdots<j_m\leq n}\lambda_{j_1}\dots\lambda_{j_m},
\end{equation}
  for $\lambda=(\lambda_1,\dots,\lambda_n)\in\mathbb{R}^n $,  $ m=1,\dots,n$.
  The function $\mathcal{H}_m$ is determined by
      \begin{equation}
  \begin{aligned}\label{BBBBB}
\det(sI+ A)&
  =\prod\limits_{k=1}^n(s+\lambda_k(A))
 =\sum\limits_{m=0}^n\sum\limits_{1\leq j_1<\dots<j_m\leq n}\lambda_{j_1}(A)\cdots\lambda_{j_m}(A)s^{n-m}=\sum\limits_{m=0}^n \mathcal{H}_m(A)s^{n-m}
   \end{aligned}
  \end{equation}for $s\in\mathbb{R}$,
   by definition.
\begin{prop}
 $\mathcal{H}_m(A)$ is a polynomial of order $m$ on $ \mathscr{H}^n$   and   is hyperbolic at $I$ for $m=1,\dots,n.$
\end{prop}
\begin{proof}
  By  Proposition $\ref{prop: hyperbolic'}$,    $\det(A)=\mathcal{H}_n$ is hyperbolic at $I$, i.e. $\det(A+sI)$ has $n$ real zeros.   If we take $Q(sI+  A)=\frac{d}{ds}\det(sI+ A)$,   the equation   $Q(sI+ A)=0$ has $(n-1)$ real zeros separating those of the equation $\det(sI + A)=0$ by Rolle's theorem \rm(cf. \cite[Lemma 1]{garding}). Thus
  \begin{equation*}
     Q(A)=Q(sI+A)\vert_{s=0}=\left. \frac{d}{ds}\right\vert_{s=0}\det(A+sI)=\mathcal{H}_{n-1}(A)
  \end{equation*}
    by   $(\ref{BBBBB})$,  and it is hyperbolic at $I$.   The result follows by repeating this procedure.
  \end{proof}
Set
\begin{equation}
 \Gamma_m : =\{A\in\mathscr{H}^n : \mathcal{H}_m(sI+A)> 0    \text{ for any  }   s\geq0\}.
\end{equation}
By definition,  $\mathcal{H}_m(sI+A)=\sum\limits_{1\leq i_1<\dots<i_m\leq n}(s+\lambda_{i_1})\cdots(s+\lambda_{i_m})> 0$ for large $s$. Then by the continuity of $\mathcal{H}_m$,   we see that
$\mathcal{H}_m(sI+A )\not= 0 $ for any $s\geq0$ if and only if $\mathcal{H}_m (sI+  A)> 0$ for any $s\geq 0$,
and so $C(\mathcal{H}_m,I)= \Gamma_m$ by definition of the cone $C(\mathcal{H}_m,I).$
\begin{cor}\label{cor:cone}
   We   have
     \begin{equation}
  \begin{aligned}
   {\Gamma}_m=\{ \mathcal{H}_1(A)>0\}\cap \cdots \{\mathcal{H}_m(A)>0\}.
   \end{aligned}
  \end{equation}
\end{cor}
\begin{proof}
It follows from  (\ref{EEEEE})  that
\begin{equation}
\mathcal{H}_m(sI+A)=\sum\limits_{1\leq i_1<\dots< i_m\leq n}(s+\lambda_{i_1})\dots(s+\lambda_{i_m})=\sum\limits_{p=0}^m\binom{n-p}{m-p}\mathcal{H}_p(A) s^{m-p}.
\end{equation}
Since $\mathcal{H}_m$ is hyperbolic at $I$, for given $A\in {\Gamma}_m$, there exist $m$ positive number $\mu_1 ,\dots,\mu_m $ such that
\begin{equation*}
\begin{aligned}
\mathcal{H}_m(sI+A)&=\binom{n}{m}\prod_{j=1}^m(s+\mu_j ) =\binom{n}{m}\sum\limits_{p=0}^m\left(\sum\limits_{1\leq i_1<\dots<i_p\leq m}\mu_{i_1} \dots\mu_{i_p} \right)s^{m-p}.
\end{aligned}
\end{equation*}  So  $\mathcal{H}_p(A)= \binom{n}{m} {\binom{n-p}{m-p}}^{-1} \sum\limits_{1\leq i_1<\dots<i_p\leq m}\mu_{i_1} \dots\mu_{i_p} >0$  for $p=1,\dots,m.$
\end{proof}
  \begin{cor}\label{garding 4}
If $A_1,\dots,A_m\in \Gamma_m$, then
\begin{equation}
\begin{aligned}
\binom{n}{m}\det(A_1,\dots,A_m,I,\dots,I)\geq \mathcal{H}_m(A_1)^{\frac{1}{m}}\dots \mathcal{H}_m(A_m)^{\frac{1}{m}}.
\end{aligned}
\end{equation}
\end{cor}
\begin{proof}   Apply Proposition  \ref{point}  to $P=\mathcal{H}_m$ to get
$$M(A_1,\dots,A_m)\geq  \mathcal{H}^{\frac{1}{m}}_m(A_1)\dots\mathcal{H}^{\frac{1}{m}}_m(A_m),$$
 where $M$ is the completely polarized form of $ \mathcal{H}_m$. Recall that the   completely polarized form $M$ of a hyperbolic polynomial $P$ is a polynomial   uniquely determined by being linear in each argument, invariant under permutations and satisfying $M(x,\dots,x)=P(x)$ \cite{garding}.
But  $\det(A_1,\dots,A_m,I,\dots,I)$ is linear in $A_1,\dots,A_m,$ and invariant under permutations, and   $\det(A,\dots,A,I,\dots,I)= {\mathcal{H}_m(A)}/{\binom{n}{m}} $ (cf. (\ref{12q})).
Therefore, \begin{equation}
M(A_1,\dots,A_m)=\binom{n }{m }\det(A_1,\dots,A_m, I,\dots, I).
\end{equation}
  The result follows.
\end{proof}

\section{The quaternionic $m$-Hessian operator }

Alesker introduced
  the   quaternionic Monge-Amp\`ere operator in \cite{alesker1}.
For a point $q=(q_0\dots q_{n-1})  \in \mathbb{H}^n$, write
 $
q_l=x_{4l}+x_{4l+1}\mathbf{i}+x_{4l+2}\mathbf{j}+x_{4l+3}\mathbf{k},
 $
$l=0,\dots,n-1.$ The {\it Cauchy-Fueter operator} is
  \begin{equation}
  \begin{aligned}
  \frac{\partial u}{\partial \overline{q_l}}=\partial_{x_{4l}}+\mathbf{i}\partial_{x_{4l+1}}+\mathbf{j}\partial_{x_{4l+2}}+\mathbf{k}\partial_{x_{4l+3}},
  \end{aligned}
  \end{equation}
  and its conjugate
   $
     \frac{\partial u}{\partial {q_l}}=\partial_{x_{4l}}-\mathbf{i}\partial_{x_{4l+1}}-\mathbf{j}\partial_{x_{4l+2}}-\mathbf{k}\partial_{x_{4l+3}}.
  $
For a $C^2$ function $u$, the {\it{quaternionic Monge-$Amp\grave{e}re$ operator}}  on $\mathbb{H}^n$ is defined as the Moore determinant of its quternionic Hessian
  \begin{equation}
  \begin{aligned}
  \det\left(\frac{\partial^2 u}{\partial \overline{q_l}\partial {q_k}} \right),
  \end{aligned}
  \end{equation}
 while the {\it quaternionic $m$-Hessian operator} $\mathcal{H}_m(u)$  is defined as
  \begin{equation}
 \mathcal{H}_m(u):=\mathcal{H}_m\left( \frac{\partial^2 u}{\partial \overline{q_l}\partial {q_k}} \right).
  \end{equation}

Let us recall that   two first-order differential operator $d_0$ and $d_1$, introduced in \cite{wan-wang},  act  on the quaternionic version of differential form. The behavior of $d_0$ and $d_1$ and $\Delta=d_0d_1$  is very similar to $\partial$,$\overline\partial$ and $\partial\overline\partial$ in several complex variables. This formulation of the   quaternionic $m$-Hessian operator is fundamental here  in the sense that we can use Stokes-type formula, etc.

  \subsection{Positive forms}
  Fix a basis $\{\omega^0,\omega^1,\dots,\omega^{2n-1}\}$ of $\mathbb{C}^{2n}$. Let $\wedge^{2k}\mathbb{C}^{2n}$ be the complex exterior algebra generated by $\mathbb{C}^{2n}$, $0\leq   k\leq   n$. Recall   the embedding
$
\tau : M_{\mathbb{H}}(p,r)\rightarrow M_{\mathbb{C}}(2p,2r)
$
as follows, where $M_{\mathbb{F}}(p,r)$ is the space of all $p\times r$-matrices over   field $\mathbb{F}$. For a quaternionic $(p\times r)$-matrice $M$,
write $\mathcal{M}=a+b\mathbf{j}$ for some complex matrices $a,b \in M_{\mathbb{C}}(p,r)$. Then
      \begin{equation}\label{eq:tau}
  \tau( {M}): =\left(
                                                           \begin{array}{cc}
                                                              a &    -b \\

                                                  \overline{b} &\overline{a}\\
                                                                                                                 \end{array} \right)
   \end{equation}
  (cf. \cite{wang2}).  We will notations in \cite{wang2}, as the relabelling of those in \cite{wan-wang}, which have advantages in the proof of some properties of quaternionic linear algebra.

For $M\in M_{\mathbb{C}}(2n,2n)$, define its {\it $\mathbb{C}$-linear action} on $\mathbb{C}^{2n}$ as  \cite{wang2}:
$
M.\omega^A=\sum\limits_{B=0}^{2n-1}M_{AB}\omega^B,
$
and   the induced action on $\wedge^{2k}\mathbb{C}^{2n}$ as
$
M.(\omega^{A_1}\wedge\dots\wedge\omega^{A_{2k}})=M.\omega^{A_1}\wedge\dots\wedge M.\omega^{A_{2k}}.
$
  For ${M}\in{M}_{\mathbb{H}}(n,n)$, defines its induced $\mathbb{C}$-linear action on $\mathbb{C}^{2n}$ as ${M}.\omega^A=\tau({M}).\omega^A,$   and so on $\wedge^{2k}\mathbb{C}^{2n}$.
     Then for ${M}\in U_{\mathbb{H}}(n)$, ${M}.\beta_n=\beta_n$ and ${M}.\Omega_{2n}=\Omega_{2n} $, where
    \begin{equation} \beta_n=\sum\limits_{l=0}^{n-1}\omega^l\wedge\omega^{n+l},\qquad \qquad \beta_n^n=\wedge^n\beta_n=n!~\Omega_{2n},
\end{equation}
where $\Omega_{2n}:=\omega^0\wedge\omega^n\dots\wedge\omega^{n-1}\wedge\omega^{2n-1} $.

  There exists a real linear action $\rho(\textbf{j})$ on $\mathbb{C}^{2n}$ \cite{wan-wang}:
\begin{equation}
\begin{aligned}
  \rho(\textbf{j}): \mathbb{C}^{2n}\rightarrow \mathbb{C}^{2n}, \qquad \rho(\textbf{j})(z\omega^k)=\overline{z}J.\omega^k, \qquad \text{
where }\qquad
J=\left(
                                                           \begin{array}{cc}
                                                                 0& I_n \\
                                                         -I_n& 0\\
                                                           \end{array} \right).
   \end{aligned}
  \end{equation}

   An element $\omega$ of $\wedge^{2k}\mathbb{C}^{2n}$ is called {\it real} if $\rho(\mathbf{{j}})\omega=\omega$. Denote by $\wedge_{\mathbb{R}}^{2k}\mathbb{C}^{2n}$ the subspace of   all real elements in $\wedge^{2k}\mathbb{C}^{2n}$, which is the counterpart of $(k,k)$-forms in complex analysis.

An element $\omega$ of $\wedge^{2n}_{\mathbb R}\mathbb{C}^{2n}$ is called {\it positive} if $\omega=\kappa\Omega_{2n}$ for some non-negative number $\kappa$. An element $\omega\in\wedge^{2k}_{\mathbb R}\mathbb{C}^{2n}$ is said to be {\it elementary strongly positive} if there exist linearly independent right $\mathbb H$-linear mappings $\eta_j:{\mathbb H}^n \to \mathbb H$, $j=1,\dots,k,$ such that
\begin{equation}
\omega=\eta_1^{*}\widetilde\omega^0\wedge\eta_1^{*}\widetilde\omega^1\wedge\cdots
\wedge\eta_k^{*}\widetilde\omega^0\wedge\eta_k^{*}\widetilde\omega^1,
\end{equation}
where $\{\widetilde\omega^0,\widetilde\omega^1 \}$ is a basis of ${\mathbb{C}}^2$ and $\eta_j^{*}:{\mathbb C}^2\to{\mathbb C}^{2n}$ is the induced $\mathbb C$-linear pulling back transformation of $\eta_j$. An element $\omega\in\wedge^{2k}_{\mathbb R}\mathbb{C}^{2n}$ is called {\it strongly positive} if it belongs to the convex cone ${SP}^{2k}{\mathbb C}^{2n}$ in $\wedge^{2k}_{\mathbb R}\mathbb{C}^{2n}$ generated by elementary strongly positive elements. An $2k$-element $\omega$ is said to be {\it positive} if for any elementary strongly positive element $\eta \in SP^{2n-2k}\mathbb{C}^{2n}$, $\omega\wedge\eta$ is positive. By definition, $ \beta_n$ is a strongly  positive
$2$-form, and
$ \beta_n^n $ is a positive
$2n$-form.

  \begin{prop}\label{prop: diagonal}  \cite[Theorem 1.1]{wang2}
  (1) For a complex skew symmetric matrix   $M =(M_{ AB})\in M_{\mathbb{C}}(2n$, $2n )$, the $2$-form
$
   \omega =\sum_{A,B=0}^{2n-1} M_{ AB}\,\omega^A\wedge\omega^B
$
   is real if and only if there exists a hyperhermitian $n\times n$-matrix $ \mathcal{{M}}=( \mathcal{M}_{jk})$,   such that
$
     M=\tau(\mathcal{M})J$.
  \\(2) When $\omega$   is   real, there exists a quaternionic unitary matrix $\mathcal{E}\in \text{U}_{\mathbb{H}}(n)$ such that
\begin{equation*}\label{eq:normal-matrix}
     \tau(\mathcal{E})^tM\tau(\mathcal{E})=\left(
      \begin{array}{cc}
                                                              0 &\mathcal{V}           \\
                                                                 - \mathcal{V}   &0
                                                           \end{array} \right),
                                                           \qquad \text {where } \qquad  \mathcal{V}=\operatorname{diag}(\nu_0,\ldots,
                                                               \nu_{n-1}),
\end{equation*} for some real numbers $\nu_0,\ldots,\nu_{n-1}$.
Namely, we can normalize $\omega$ as
$
    \omega =2\sum_{ l=0}^{n-1} \nu_{ l}\widetilde{\omega}^l\wedge\widetilde{\omega}^{l+n}
 $    with $\widetilde {\omega}^A= \mathcal{E^*}.{\omega}^A$. In particular, $\omega$ is (strongly) positive if and only if each $\nu_{ l}\geq 0$ ($>0$).
\end{prop}

\begin{prop}\label{prop:positive form}  \cite[Lemma 3.3]{wan-wang1}
For $\eta\in\wedge_{\mathbb{R}}^{2k}\mathbb{C}^{2n}$ with $\|\eta\|\leq 1$, $\beta_n^k\pm \epsilon\eta$ is positive $2$k-form for some sufficiently small absolute constant $\epsilon>0$.
\end{prop}

   \subsection{$d_0,d_1$ formulation of the quaternionic $m$-Hessian operator}
 We express the quaternionic $m$-Hessian operator in terms of $d_1,d_1.$   Let $\Omega$ be a domain in $\mathbb{H}^n$. Denoted by ${\mathcal{D}}^{p}({\Omega})$ the set of all $C_0^{\infty}(\Omega)$
functions valued in $\wedge^{p}\mathbb{C}^{2n}$. $F\in{\mathcal{D}}^{2k}({\Omega})$ is called a {\it (strongly) positive form} if for any $q\in\Omega$, $F(q)$ is a (strongly) positive element.     Define
  $d_0,d_1:C ^{1}(\Omega,\wedge^p\mathbb{C}^{2n})\to C (\Omega,\wedge^{p+1}\mathbb{C}^{2n})$   by
\begin{equation}
\begin{aligned}
  d_\alpha F=\sum_{I}\sum_{A=0}^{2n-1} \nabla_{A\alpha}f_I \omega^A\wedge\omega^I,
   \end{aligned}
  \end{equation}
for   $F=\sum_I f_I\omega^I \in C ^{1}(\Omega,\wedge^p\mathcal{C}^{2n})$, where the multi-index $I=(i_1\dots i_p)$, $\omega^I=\omega^{i_1}\wedge\dots \wedge \omega^{i_p}$, and
  the first-order differential operators $\nabla_{A\alpha}$ ($A=0,\ldots, 2n-1,$ $\alpha =0 ,1 $) are
\begin{equation}\label{nabla3}\left(
                                                           \begin{array}{cc}
                                                               \nabla_{00} &      \nabla_{01} \\
                                                         \vdots& \vdots\\
                                                         \nabla_{l0} &\nabla_{l1}\\
                                                         \vdots& \vdots\\
                                                           \nabla_{n0} &      \nabla_{n1} \\
                                                         \vdots& \vdots\\
                                                         \nabla_{(n+l)0} &\nabla_{(n+l)1}\\
                                                         \vdots& \vdots
                                                           \end{array} \right)=\left(
                                                           \begin{array}{cc}
                                                               \partial_{x_{0}}+\mathbf{i} \partial_{x_{1}} &   - \partial_{x_{2}}-\mathbf{i} \partial_{x_{3}} \\
                                                         \vdots& \vdots\\
                                                      \partial_{x_{4l}}+\mathbf{i} \partial_{x_{4l+1}} &-\partial_{x_{4l+2}}-\mathbf{i} \partial_{x_{4l+3}}\\
                                                         \vdots& \vdots\\
                                                         \partial_{x_{2}}-\mathbf{i} \partial_{x_{3}} &    \partial_{x_{0}}-\mathbf{i} \partial_{x_{1}}\\
                                                         \vdots& \vdots\\
                                                   \partial_{x_{4l+2}}-\mathbf{i} \partial_{x_{4l+3}}&\partial_{x_{4l}}-\mathbf{i} \partial_{x_{4l+1}}\\
                                                         \vdots& \vdots
                                                           \end{array} \right).
\end{equation}
\begin{prop} \label{prop:d^2-identities} \cite[Proposition 2.2]{wan-wang}
  (1) $d_0d_1=-d_1d_0$; \\
 (2) $ d_0 ^2= d_1 ^2=0$; \\
  (3) For $F\in C ^1(\Omega,\wedge^p\mathbb{C}^{2n}),G\in C ^1(\Omega,\wedge^q\mathbb{C}^{2n})$,we have
   \begin{equation*}
 d_\alpha(F\wedge G)=d_\alpha F\wedge G+(-1)^pF\wedge d_\alpha G,\qquad \alpha=0,1.
   \end{equation*}
  \end{prop}
The following nice identity will be frequently used.
  \begin{prop} \label{prop:dd-identities}  \cite[Proposition 2.3]{wan-wang}
  For $u_1,\dots,u_n \in C^2$,
\begin{equation}
\begin{aligned}
\Delta u_1\wedge\Delta u_2\wedge \dots\wedge\Delta u_n
&=d_0(d_1u_1\wedge\Delta u_2\wedge \dots\wedge\Delta u_n)=-d_1(d_0u_1\wedge\Delta u_2\wedge \dots\wedge\Delta u_n)\\
&=d_0d_1(u_1\Delta u_2\wedge \dots\wedge\Delta u_n)
=\Delta(u_1\Delta u_2\wedge \dots\wedge\Delta u_n).\label{C}
\end{aligned}
\end{equation}
\end{prop}

 Define
   \begin{equation*}
     \int_\Omega F=\int_\Omega fdV,
  \end{equation*}
if   $F=f\Omega_{2n}\in L^1(\Omega,\wedge^{2n} \mathbb{C}^{2n})$, where $dV$ is the Lebesgue measure.

\begin{lem}\label{lem:Stokes} \cite[Lemma 3.2]{wan-wang}  {\rm (Stokes-type formula)}  Assume that $T=\sum_AT_A\omega^{\widehat{A}}$ is a $C^1$ $(2n-1)$-form in $\Omega$, where $\omega^{\widehat{A}}=\omega^A  \rfloor \Omega_{2n}:=(-1)^{A-1}\omega^0\wedge\ldots\wedge\omega^{A-1}\wedge\omega^{A+1}\wedge\ldots\wedge\omega^{2n-1}$. Then for a $C^1$ function $h$, we have
\begin{equation}\label{eq:stokes}
\int_\Omega hd_\alpha T=-\int_\Omega d_\alpha h\wedge T+\int_{\partial\Omega}\sum_{A=0}^{2n-1}h\,T_A\, \tau( \mathbf{{n}})_{A\alpha}\,
dS,\qquad \alpha=0,1,
\end{equation}
where
 $  \mathbf{{n}}: =(n_0,n_1,\ldots,n_{4n-1})$ is the unit outer normal vector to $\partial\Omega$, $dS$ denotes the surface measure of $\partial\Omega$, and  $\tau( \mathbf{{n}})$ is a complex $(2n)\times 2$-matrix by definition \eqref{eq:tau} of $\tau$.  In particular, if $h=0$ on $\partial\Omega$, \eqref{eq:stokes} has no boundary term.
\end{lem}
Recall   the {\it Baston operator}
$
  \Delta u:=d_0d_1 u$ for a real $C^2$ function $u$.

\begin{prop}\label{relationship}   \cite[Theorem 1.3]{wan-wang}
Let $u_1,\dots,u_n $ be real $C^2$ functions on ${\mathbb{H}}^n$. Then we have
\begin{equation}
\Delta u_1\wedge \dots\wedge\Delta u_n=n!\det(A_1,A_2,\dots,A_n)\Omega_{2n}.
\end{equation}
where  $A_j=\Big(\frac{\partial^2 u_j}{\partial \overline{q_l}\partial {q_k}}(q)\Big) $.
\end{prop}

\begin{prop} \label{hessian operator}    \begin{equation}
  (\Delta u)^{m}\wedge  \beta_n ^{n-m}
=m!(n-m)!\mathcal{H}_m(u)\Omega_{2n}.
  \end{equation}
  \end{prop}
\begin{proof}
  Apply Proposition $\ref{relationship}$ to $u_1= \cdots=u_m=u$ and $u_{m+1}=\cdots=u_n=\|q\|^2$ to get
  \begin{equation}
  8^{n-m} (\Delta u)^m\wedge  \beta_n ^{n-m} =n!\det (A,\dots,A,8I,\dots,8I )\Omega_{2n},
  \end{equation}
  where
   $A=\Big(\frac{\partial^2 u}{\partial \overline{q_l}\partial {q_k}}(q)\Big),$ and
   \begin{equation}
      \Delta \|q\|^2=d_0d_1\|q\|^2=8{\beta }_n.
   \end{equation}
  By definition,   the coefficient of the monomial $s_1\dots s_n$ of $\det (s_1A+\dots+s_mA+8s_{m+1}I+\dots +8s_n I )$   divided by $n!$ is the $\det (A,\dots,A,8I,\dots,8I )$.
On the other hand,   we can find a   quaternionic unitary matrix   $\mathcal{U}  \in \text{U}_{\mathbb{H}}(n)$ such that
$
{ { \mathcal{U}}}^*A\mathcal{U}=\operatorname{diag}(\lambda_1,\dots,\lambda_n).
$ Now apply Proposition  \ref{diagonal and real} to get
  \begin{equation}
  \begin{aligned}\label{12q}
  \det \left (\sum_{j=1}^m s_j A+8\sum_{j=  m+1}^n s_j I \right)&
=\det\left(\mathcal{U}^*\left (\sum_{j=1}^m s_j A+8\sum_{j=  m+1}^n s_jI\right)\mathcal{U}\right)\\
  &=\det\left (\sum_{j=1}^m s_j\operatorname{diag}(\lambda_1,\dots,\lambda_n)  +8\sum_{j=  m+1}^n s_j I\right)\\
       &=\prod\limits_{p=1}^n\left (\lambda_p\sum_{j=1}^m s_j+8\sum_{j=  m+1}^n s_j\right),                                                                                                                       \end{aligned}
  \end{equation}
whose coefficient of ${s_1}\dots{s_n}$ is
$
8^{n-m}m!(n-m)!\sum_{1\leq i_{1} \leq \dots \leq i_{m} \leq n}\lambda_{i_1}\dots \lambda_{i_m}.
$
Therefore
   \begin{equation}
  \begin{aligned}\label{kkkk}
  (\Delta u)^{m}\wedge  \beta_n ^{n-m}
=m!(n-m)!\sum_{1\leq i_1 \leq \dots \leq i_m \leq n}\lambda_{i_1}\dots \lambda_{i_m}\Omega_{2n}.
\end{aligned}
  \end{equation}
The result follows.
  \end{proof}
We also need the following elementary strong  positivity (cf., e.g. \cite[Proposition 4.2]{wang21}).
\begin{prop} \label{prop:wedge-positive}   For any $C^1$ real function $u$,  $d_0u\wedge d_1u$ is elementary strongly positive if grad $u\neq0$.
\end{prop}

\section{\textbf{Quaternionic $m$-subharmonic functions}}
\subsection{Smooth  quaternionic $m$-subharmonic function}
A  real   $C^2$ functions $u$   is said to be {\it{quaternionic $m$-subharmonic}} on    $\Omega\subset\mathbb{H}^n$   if
\begin{equation}\label{eq:QSH-positivity}
     \left(\frac{\partial^2 u}{\partial \overline{q_l}\partial {q_k}}\right)(q)\in \overline{\Gamma}_m
  \end{equation}
 for any $q\in \Omega$.   It follow from Corollary \ref{cor:cone} and Proposition $\ref{hessian operator}$ that it is equivalent to require \begin{equation}
  \begin{aligned}
 (\Delta u)^{k}\wedge  \beta_n ^{n-k}\geq 0,    \qquad \text{  for } k=1,2, \dots,m.
\end{aligned}
  \end{equation}

   \begin{prop}\label{prop:QSH-positivity}
If $u_1,\dots,u_k$ are $ C^2 $  quaternionic $m$-subharmonic functions, $1\leq k\leq m$, then $\Delta u_1\wedge\dots\wedge\Delta u_k\wedge\beta_n^{n-m}\geq  0.$
\end{prop}
\begin{proof} Since $u_1,\dots,u_m \in QSH_m(\Omega)\cap C^2(\Omega)$, $A_1=\left(\frac{\partial^2 u_1}{\partial \overline{q_l}\partial {q_k}}\right),\dots,A_m=\left(\frac{\partial^2 u_m}{\partial \overline{q_l}\partial {q_k}}\right)\in \overline{\Gamma}_m$.   Then we have
   \begin{equation*}
  \begin{aligned}
\binom{n}{m}\det(A_1,\dots,A_m,I,\dots,I)\geq \mathcal{H}_m(A_1)^{\frac{1}{m}}\dots\mathcal{H}_m(A_m)^{\frac{1}{m}}\geq0
\end{aligned}
  \end{equation*}
by Garding's inequality in Corollary $\ref{garding 4}$. Then,  by Proposition $\ref{relationship}$, we get
     \begin{equation*}\label{eq:positivy-m}
  \begin{aligned}
\Delta u_1\wedge\dots\wedge\Delta u_m\wedge\beta_n^{n-m}=n!\det(A_1,\dots,A_m,I,\dots,I)\Omega_{2n}\geq0.
\end{aligned}
  \end{equation*}

  For $k<m$, it is sufficient to prove that
\begin{equation}
\begin{aligned}\label{another equality}
\Delta u_1\wedge\dots\wedge\Delta u_k\wedge\beta_n^{n-m}\wedge\omega\geq  0.
\end{aligned}
\end{equation}for any elementary strongly  positive $2(m-k)$-element $ \omega=\eta^*_1\widetilde{\omega}^0\wedge\eta^*_1\widetilde{\omega}^1
\wedge\dots\wedge\eta^*_{m-k}\widetilde{\omega}^0\wedge\eta^*_{m-k}\widetilde{\omega}^1, $ where $\eta_j:\mathbb{H}^n\to\mathbb{H},j=1,\dots,m-k,$ are linearly independent right $\mathbb{H}$-linear mappings and $\{\widetilde{\omega}^0,\widetilde{\omega}^1\}$ is a basis of $\mathbb{C}^2.$
Since $\Delta \|\widetilde{q}_0\|^2  =8\widetilde{\omega}^0\wedge\widetilde{\omega}^1$ and $\eta^*_j(\Delta \|\widetilde{q}_0\|^2  )=\Delta(  {\| \eta_j ( q)\|^2}  )$. So $(\ref{another equality})$ is proved by $\eta_j(  q  )\in QPSH\subset QSH_m(\Omega)$ and the case $k=m$ in  \eqref{eq:positivy-m}.
\end{proof}
\subsection{Closed positve currents}\rm
To define nonsmooth quaternionic $m$-subharmonic functions, we need to use   currents.   An element of the dual space ($ \mathcal{D} ^{2n-p}({\Omega})'$) is called a {\it $p$-current}.   Obviously   $2n$-currents are just   distributions on $\Omega$. A $2k$-current $T$ is said to be {\it positive} if we have $T(\eta)\geq 0$ for any strongly positive form $\eta\in{\mathcal{D}}^{2n-2k}({\Omega})$. Let $\psi$ be a $p$-form whose coefficients are locally integrable in $\Omega$. One can associate with $\psi$ the $p$-current $T_\psi$ defined by
 $
T_\psi(\varphi)=\int_\Omega \psi\wedge\varphi
$
  for   $\varphi\in{\mathcal{D}}^{2n-p}({\Omega})$.

Now for a $p$-current $F$, we define the $(p+1)$-current $d_\alpha F$ as
\begin{equation*}\label{eq:generalized-sense}(d_\alpha F)(\eta):=-F(d_\alpha\eta),\qquad \alpha=0,1,
\end{equation*}for any test form $\eta\in \mathcal{D} ^{2n-p-1}({\Omega})$. We say a form (or a current) $F$ is \emph{closed} if
$d_0F=d_1F=0 .$

If a $p$-current $T$ has a continuous extension to the space of $ (2n-p)$-forms
with continuous coefficients,   it is called a $p$-current of {\it order zero} or
{\it of measure type}. A $p$-current $T$ is of measure type if and only if
for any  neighborhood $G\Subset \Omega$, there
exists a constant $K_G $ such that
$
   |T(\alpha)|\leq K_G \|\alpha\|_G,
$
where $\|\alpha\|_G=\sum_{ I }' \max_{ q\in G}|\alpha_I(q)|$ for $\alpha=\sum_{|I|=2n-p}' \alpha_I\omega^I$. Here the summation $\sum '$ is taken over increasing indices of length $2n-p$.

Denote by $\mathcal{M}^p(\Omega)$ the set of all $p$-currents of measure type, and it is identified with
     $\wedge^p $-valued Radon measures on $\Omega$. A sequence of currents
 $ T_j\in \mathcal{M}^p(\Omega)$ {\it weakly $*$ converges} to  $T$ if
 $T_j (\alpha)\rightarrow T(\alpha)$
for any $ (2n-p)$-forms
with continuous coefficients.  A family of currents
 $ T_\kappa\in \mathcal{M}^p(\Omega)$ is {\it weakly $*$ compact} (or {\it locally uniformly bounded}) if and
only if for any domain $G \Subset \Omega$ there is a constant $K_G $ depending only on $ G $
such that
\begin{equation}
   |T_\kappa(\alpha)|\leq K_G \|\alpha\|_G.
\end{equation}
\subsection{ Non-smooth quaternionic $m$-subharmonic functions}
A   $[-\infty,\infty)$-valued upper semicontinuous   function $u\in{L}_{loc}^1(\Omega)$ is called {\it{quaternionic $m$-subharmonic}}, if for any $   C^4  $ quaternionic $m$-subharmonic functions $v_1,\dots,v_{m-1}$ on $\Omega$, the current   $\Delta u\wedge\Delta v_1\wedge\dots\wedge \Delta v_{m-1}\wedge\beta_n^{n-m}$ defined by
\begin{equation}\label{eq:Qmsubharmonic}
\Delta u\wedge\Delta v_1\wedge\dots\wedge \Delta v_{m-1}\wedge\beta_n^{n-m}(\omega) =\int u\Delta v_1\wedge\dots\wedge \Delta v_{m-1}\wedge\beta_n^{n-m}\wedge\Delta \omega,\qquad \text{ for any} \quad \omega \in C_0^{\infty}(\Omega),
  \end{equation}
  is nonnegative. The set of quaternionic $m$-subharmonic functions on $\Omega$ is denoted by $QSH_m(\Omega)$.

 \begin{prop}
A function $u\in C^2(\Omega)$ is   quaternionic $m$-subharmonic in the above sense  if and only if \eqref{eq:QSH-positivity} holds for any $q\in \Omega$.
   \end{prop}
  \begin{proof} For a function $u\in C^4(\Omega)$,
  \begin{equation}\label{eq:Qmsubharmonic2}
     \int_\Omega u\Delta v_1\wedge\dots\wedge \Delta v_{m-1}\wedge\beta_n^{n-m}\wedge\Delta \omega=\int_\Omega \omega\Delta u\wedge\Delta v_1\wedge\dots\wedge \Delta v_{m-1}\wedge\beta_n^{n-m}
  \end{equation}by applying Stokes-type formula \eqref{eq:stokes} twice, since integrands  vanish on the boundary.
  By continuity,
\eqref{eq:Qmsubharmonic2}   is nonnegative for any nonnegative $\omega$  if and only if  $\Delta  u \wedge\Delta v_1\wedge\dots\wedge \Delta v_{m-1}\wedge\beta_{n }^{n-m}$ is positive at  each $q \in \Omega$. So in this case, the definition \eqref{eq:Qmsubharmonic}  is equivalent  to require $v_1,\dots,v_{m-1}$ only to be quadratic $ QSH_m  $ polynomials.

 {\it Sufficiency}. By Proposition  \ref{prop:QSH-positivity}, $\Delta u\wedge\Delta v_1\wedge\dots\wedge \Delta v_{m-1}\wedge\beta_n^{n-m} $ in \eqref{eq:Qmsubharmonic}  is a positive form if the positivity in \eqref{eq:QSH-positivity} holds for $u$.

  {\it Necessity}. We prove it by induction on dimension $n$ of the space and the number $m$. Suppose that we have proved the result for dimension less than $n$ and $m-1$ on dimension   $n$.
    Now by rotation if necessary, we can assume that $\left(\frac{\partial^2 u}{\partial \overline{q_l}\partial {q_k}}\right)(q_0)$ is diagonalized with eigenvalues $\lambda_1(q_0)\leq\cdots\leq \lambda_n(q_0)$. Hence $\lambda_n(q_0)\geq 0$ and
  \begin{equation}\label{eq:Sm-eigenvalue}
   \mathcal{H}_m(u)(q_0) = \lambda_{n}(q_0) \sum\limits_{1\leq j_2<\dots<j_m\leq n-1}\lambda_{j_1}(q_0) \cdots\lambda_{j_{m-1}}(q_0)+ \sum\limits_{1\leq j_1<\dots<j_m\leq n-1}\lambda_{j_1}(q_0)\cdots\lambda_{j_m}(q_0).
  \end{equation}
   If we take $\triangle v_{m-1}=\omega^{n-1}\wedge \omega^{2n-1}$, i.e. $v_{m-1}=|q_n|^2$. Then the positivity of $\Delta u\wedge\Delta v_1\wedge\dots\wedge \Delta v_{m-1}\wedge\beta_n^{n-m} $ at point $q_0$ implies that
   \begin{equation*}
      \Delta' u(q_0 )\wedge\Delta' v_1(q_0 )\wedge\dots\wedge \Delta' v_{m-2}(q_0 )\wedge\beta_{n-1}^{n-m}
   \end{equation*}
    is a positive  element on $\mathbb{H}^{n-1}$, where $\Delta'$ is the Baston operator on $\mathbb{H}^{n-1}$.
By the assumption of induction for dimension $n-1$, we see that $\left(\frac{\partial^2 u}{\partial \overline{q_l}\partial {q_k}}(q_0)\right)_{1\leq j,k\leq n-1}$ belongs to $  \overline{\Gamma}_{m-1}$. Thus, the second sum  in \eqref{eq:Sm-eigenvalue} is non negative. The first sum  in \eqref{eq:Sm-eigenvalue} is also non negative by the   assumption of induction for $m-1$ in dimension $n$.
\end{proof}

  \begin{prop} \label{prop:QSH-m}
Let $\Omega$ be a domain in $\mathbb{H}^n$.  Then,
\\
(1)  The   standard approximation   $u_\epsilon=u\ast\chi_\epsilon$ is also a  $QSH_m$ function, and satisfies $u_\epsilon\downarrow u$ as  $\epsilon\downarrow 0$.\\
(2) $QPSH=QSH_n\subset\dots\subset QSH_1=SH$.\\
(3) $au+bv \in QSH_m(\Omega)$  for any $a,b\geq  0$.\\
(4) If $\gamma(t)$ is a convex increasing function   on $\mathbb{R}$ and $u\in QSH_m$, then $\gamma\circ u \in QSH_m$. \\
(5) The limit of a uniformly converging or decreasing sequence of $QSH_m$ functions is an $QSH_m$ function.\\
(6) The maximum of a finite number of $QSH_m$ functions is a  $QSH_m$ function; for an arbitrary locally uniformly bounded family $\{u_{\alpha}\}\subset QSH_m$, the regularization $u^*(q)$ of the supremum $u(q) = \sup_{\alpha}u_{\alpha}(q)$ is also a  $QSH_m$ function.
\\(7)
If $D$ is an open subset  of $\Omega$, $ u \in  QSH_m(\Omega)$, $ v \in  QSH_m  (D)$  and $\limsup_{q\rightarrow q_0} v(q ) \leq u(q_0)$ for all $q_0\in \partial D\cap \Omega$,
then the function defined by
 \begin{equation}\label{eq:paste}
   \phi=\left\{
   \begin{array}{ll}
u, \qquad &{\rm on}\quad \Omega \setminus D,\\ \max\{u, v\}, \quad &{\rm on}\quad  D,
\end{array}
\right.
 \end{equation}
belongs to $QSH_m(\Omega)$.
\end{prop}
\begin{proof} Because there is no characterization of $m$-subharmonicity by the  submean value inequality, the proof is different from that for plurisubharmonic functions.

(1)  For any $ C^4(\Omega)\cap QSH_m(\Omega)$ functions $v_1,\dots,v_{m-1}$ and nonnegative   function $\omega\in C_0^\infty(\Omega)$, it is direct to see that if $\epsilon>0$ small,
\begin{equation}
  \begin{aligned}\label{YYYYY}
&\int_{\Omega} \Delta u_{\epsilon}(x)\wedge\Delta v_1(x)\wedge\dots\wedge\Delta v_{m-1}(x)\wedge \beta_n^{n-m}\wedge \omega(x) \\
=&\int_{B(0,\epsilon)}\chi_{\epsilon}(y)dV(y)  \int_{\Omega} u(z)\Delta v_1( z+y)\wedge\dots\wedge\Delta v_{m-1}(z+y )\wedge\beta_n^{n-m}\wedge\Delta\omega(z+y ) \geq 0,
\end{aligned}
  \end{equation} by \eqref{eq:Qmsubharmonic} for $u$ with $\omega(\cdot  ) $ replaced by  $\omega(\cdot+y ) $ and $v_{j} $ replaced by  $ v_{j}(\cdot +y ) $. Thus $ u_{\epsilon}$ is $ QSH_m   $.

For  $v_1,\dots,v_{m-1}\in C^4(\Omega)\cap QSH_m(\Omega)$, denote
\begin{equation}\label{eq:alpha}
   \alpha := \Delta v_1\wedge\dots\wedge \Delta v_{m-1}\wedge\beta_n^{n-m}.
\end{equation}Then the linear operator $ \mathcal{A}_\alpha $
defined by
\begin{equation*}\label{eq:A-alpha}
   \mathcal{A}_\alpha( u)\cdot  \Omega_{2n}= \Delta u\wedge\Delta v_1\wedge\dots\wedge \Delta v_{m-1}\wedge\beta_n^{n-m}
\end{equation*}
is a differential operator of the second order with $C^2$ coefficients, whose symbol $\sigma(  \mathcal{A}_\alpha)(\xi; q)$ at point $q$ and direction $0\neq\xi\in \mathbb{R}^{4n}$ is given by
\begin{equation*}
 \sigma( \mathcal{A_\alpha})(\xi; q)  \Omega_{2n}= d_0|\xi|^2\wedge  d_1|\xi|^2\wedge\omega_1 \wedge\dots\wedge \omega_{m-1}\wedge\beta_n^{n-m}\geq 0
\end{equation*} where $\omega_j=\Delta v_j(q)$,
 and $d_0|\xi|^2\wedge  d_1|\xi|^2$ is elementary strongly positive by Proposition \ref{prop:wedge-positive}. Without loss of generality, we may assume the it is strictly positive, i.e. $ \mathcal{A}_\alpha $ is a uniform elliptic operator.  Otherwise, we replaced $v_j(q)$ by $v_j(q)+\varepsilon|q|^2$.
  It is also an operator of divergence form, which can be proved by  $\mathcal{A}_\alpha( u)\cdot  \Omega_{2n}= d_0(d_1 u\wedge\Delta v_1\wedge\dots\wedge \Delta v_{m-1}\wedge\beta_n^{n-m})$ by  Proposition \ref{prop:dd-identities}.

Now the positivity of \eqref{eq:Qmsubharmonic} is equivalent to   $ \mathcal{A}_\alpha u \geq 0$ in the sense of distributions, i.e. $u $ is $ \mathcal{A}_\alpha$-subharmonic. It is well known $ \mathcal{A}_\alpha$-subharmonicity can be characterized as the maximum principle, i.e.
 for every domain $G\Subset \Omega$, if $v \in C(\overline{G} ) $ satisfies $\mathcal{A}_\alpha v=0$ and
$ u \leq v  $ on $\partial G$, then $ u \leq v $ in $G$.

  All other properties can be proved by using this characterization and well known corresponding  properties for $ \mathcal{A}_\alpha$-subharmonic functions (cf. e.g. \cite{HKM}), since $\mathcal{A}_\alpha $ is an elliptic  differential operator of the second order with $C^2$ coefficients and  of divergence form.

  For example, for $ u \in  QSH_m(\Omega)$ and $ v \in  QSH_m  (D)$, they are $ \mathcal{A}_\alpha$-subharmonic on $\Omega$ and $D$, respectively. Then the function $\phi$ in \eqref{eq:paste}
   is also $ \mathcal{A}_\alpha$-subharmonic on $\Omega$ for any $\alpha := \Delta v_1\wedge\dots\wedge \Delta v_{m-1}\wedge\beta_n^{n-m}$ with $v_1,\dots,v_{m-1}\in C^4(\Omega)\cap QSH_m(\Omega)$. Thus \eqref{eq:Qmsubharmonic2}   is nonnegative for any nonnegative $\omega$.

  If $ \mathcal{A}_\alpha $ is not uniformly elliptic, we use $ \mathcal{A}_{\alpha_\epsilon} $, where $\alpha_\epsilon$ is the  $\alpha $  in \eqref{eq:alpha} with \ $v_j(q)$ replaced by $v_j(q)+\varepsilon|q|^2$. Since $ \mathcal{A}_{\alpha_\epsilon} $ is   uniformly elliptic,    $ \mathcal{A}_{\alpha_\epsilon}\phi\geq 0 $ in the sense of distributions. Then  $ \mathcal{A}_{\alpha }\phi\geq 0 $  by letting $\varepsilon\rightarrow0$. Thus, $\phi$ belongs to $QSH_m(\Omega)$ by definition.
\end{proof}

\begin{rem} (1) In the definition of $QSH_m $,  we require  $v_{j}\in C^4 $ instead of the usual condition $v_{j}\in C^2 $  in order to
 make  $ \mathcal{A}_\alpha$ of  $C^2$ coefficients.

 (2) In the complex case, the proof of these properties were only sketched in \cite{sadullaev1}, as far as I know, by using integral representation formula of solutions to the operator $ \mathcal{A}_\alpha $. But there is also the degenerate problem there.
\end{rem}

A set $E\subset \Omega$ is said to be {\it quaternionic $m$-polar} in $\Omega$, if there exists a function $u \in QSH_m(\Omega)$  such that $u \not \equiv -\infty$ and $u\vert_{E}\equiv-\infty$.
\section{Quaternionic $m$-Hessian measure and the  comparison principle}
\subsection{Quaternionic $m$-Hessian measure}
We need the following coaera formula.
\begin{prop}   \cite[Theorem 1.2.4]{Ma}
   For
     a    measurable nonnegative function $\Phi$ on an open subset $\Omega$
  of
$  \mathbb{R} ^N$   and
  $f\in C^{0,1}(\Omega)$, we have
\begin{equation} \label{eq:gradient} \int_\Omega \Phi(x)\,|\operatorname{grad} f(x)|\,dV(x)
=\int_0^\infty ds\int_{  \Omega\cap
\{|f|= s \}}\Phi(x)\,dS(x),
\end{equation} where $dS$ is the $(N-1)$-dimension Hausdorff measure $d\mathcal{H}^{N-1}$, which equals to the surface measure if the surface is smooth.
\end{prop}

A domain $\Omega$ is called {\it   $m$-hyperconvex}  if there exists
a continuous function  $\varrho\in    QSH_m(\Omega)$  such that $\varrho < 0 $ in $\Omega$ and $\lim_{q\rightarrow \partial\Omega} \varrho(q)=0$, i.e. $\{\varrho(q)<c\}$ is relatively compact in $ \Omega$ for any $c < 0$. It is called  {\it strongly $m$-hyperconvex}
if  $\varrho\in    QSH_m(G)$ for some open set $  G \Supset\Omega$.
We need the following key integral estimate. See
Sadullaev-Abdullaev
\cite[Theorem 16.2]{sadullaev0} for plurisubharmonic   functions and \cite{sadullaev1} for $m$-subharmonic functions on a ball.

\begin{thm} \label{thm:coarea} Let $  \Omega=\{\varrho<0\} $ be  a   $m$-hyperconvex domain with $\varrho\in C^2(\Omega)$,
 $\sigma=\min_\Omega \varrho$. For $u_1 \cdots  u_k \in QSH_m(\Omega)\cap C(\Omega)$, $k=0,\ldots, m$, and any $\sigma<r<0$,
\begin{equation}\label{eq:coarea}
\int_\sigma^rdt  \int_{\varrho \leq t}(\Delta \varrho)^{n-k}\wedge \Delta u_1\wedge\dots\wedge \Delta u_k \leq (M-M') \int_{\varrho \leq r}(\Delta \varrho)^{n-k+1}\wedge \Delta u_1\wedge\dots\wedge \Delta u_  { k-1 } ,
\end{equation}
where $M=\max_{\varrho \leq r}\{u_1 ,\ldots u_k \}$, $M'=\min_{\varrho \leq r} \{u_1 ,\ldots u_k \}$. In particular, if $u_k|_{\varrho = r }=0$, we have
\begin{equation}\label{eq:coarea2}
\int_\sigma^rdt  \int_{\varrho \leq t}(\Delta \varrho)^{n-k}\wedge \Delta u_1\wedge\dots\wedge \Delta u_k =- \int_{\varrho \leq r}u_k(\Delta \varrho)^{n-k+1}\wedge \Delta u_1\wedge\dots\wedge \Delta u_  { k-1 } .
\end{equation}
\end{thm}
We first prove the result under the $  C^2 $ assumption.
\begin{lem} \label{lem:coarea} Theorem \ref{thm:coarea} holds for
   $u\in QSH_m(B)\cap C^2(\Omega)$.
\end{lem}
\begin{proof} Note that  $  \mathbf {n}  = \operatorname{grad} \varrho /{|\operatorname{grad}\varrho|}$ and so $\tau( \mathbf{{n}})_{A\alpha}= {\nabla_{A\alpha}\varrho}/{|\operatorname{grad}\varrho|}$. Denote $\Theta:=\Delta u_1\wedge\dots\wedge \Delta u_ { k-1 } $.  Apply Proposition \ref{prop:dd-identities},  Stokes-type formula \eqref{eq:stokes}   and the coaera formula
\eqref{eq:gradient} to get
   \begin{equation*}\begin{split}
   \int_\sigma^rdt  \int_{\varrho \leq t}(\Delta \varrho)^{n-k}\wedge \Delta u_1\wedge\dots\wedge \Delta u_k &=\int_\sigma^rdt  \int_{\varrho \leq t}d_0\left(d_1u_k\wedge(\Delta \varrho)^{n-k}\wedge (\Delta u)^{k-1} \right)\\ &=\int_\sigma^rdt  \int_{\varrho = t}\sum_{A=0}^{2n-1}  \left(d_1u_k\wedge \Theta \wedge(\Delta \varrho)^{n-k}\right)_A\frac {\nabla_{A0}\varrho \,dS}{|\operatorname{grad}\varrho|}\\
  & =  \int_{\varrho \leq  r}\sum_{A=0}^{2n-1}  (d_1u_k\wedge \Theta\wedge(\Delta \varrho)^{n-k})_A\nabla_{A0}\varrho\, dV\\
  &=-\int_{\varrho \leq  r} d_1u_k\wedge d_0\varrho\wedge  \Theta \wedge(\Delta \varrho)^{n-k}
  \\
  &=-\int_{\varrho =  r} u_k\sum_{A=0}^{2n-1} \left(d_0\varrho\wedge \Theta \wedge(\Delta \varrho)^{n-k}\right)_A\tau( \mathbf{{n}})_{A 1} dS\\
  &\quad-\int_{\varrho \leq  r}  u_k \Theta \wedge(\Delta \varrho)^{n-k+1}:=I_1+I_2.
\end{split}
\end{equation*}In the forth identity, we have used
\begin{equation} \label{eq:boundary-term}
   \sum_{A=0}^{2n-1} \nabla_{A\alpha }\varrho\left(d_1u_k\wedge \Theta\wedge(\Delta \varrho)^{n-k}\right)_A \Omega_ {2n }=d_\alpha\varrho\wedge d_1u_k\wedge \Theta \wedge(\Delta \varrho)^{n-k},
\end{equation}
since $d_\alpha\varrho=\sum_{A=0}^{2n-1}\nabla_{A\alpha}\varrho \,\omega^A$.
 But
 \begin{equation*}
   -\sum_{A=0}^{2n-1} \tau( \mathbf{{n}})_{A 1} \left(d_0\varrho\wedge \Theta \wedge(\Delta \varrho)^{n-k}\right)_A \Omega_ {2n }=d_0\varrho\wedge d_1 \varrho\wedge \Delta u_1\wedge\dots\wedge \Delta u_  { k-1 } \wedge(\Delta \varrho)^{n-k}/{|\operatorname{grad}\varrho|}
\end{equation*}
 is nonnegative by
  using Proposition \ref{prop:wedge-positive} and \ref{prop:QSH-positivity}. So we have
\begin{equation*}\begin{split}
 I_1   &\leq - M\int_{\varrho =  r} \sum_{A=0}^{2n-1}   \left(d_0\varrho\wedge  \Theta \wedge(\Delta \varrho)^{n-k}\right)_A\tau( \mathbf{{n}})_{A 1} dS
 =M \int_{\varrho \leq  r}    \Theta\wedge(\Delta \varrho)^{n-k},
\end{split}
\end{equation*}
and
\begin{equation*}
   I_2\leq -M'\int_{\varrho \leq  r}  \Theta \wedge(\Delta \varrho)^{n-k+1}.
\end{equation*}The estimate follows. If $u_k|_{\varrho = r }=0$, we get $I_1=0$.
\end{proof}
Applying   \eqref{eq:coarea} to the ball $B=B(0,1)$ with $\varrho(q)=|q|^2-1$ repeatedly, we get
   \begin{equation*}\label{eq:iit-Hessian}\begin{split}&
\int_0^1dt_1  \int_0^{t_1}dt_2\cdots\int_0^{t_{k-1}}dt_k\int_{|q|^2\leq t_k}\Delta u_1\wedge\dots\wedge \Delta u_k\wedge\beta_n^{n-k}
\leq  (M-M')^k \int_{|q|^2\leq 1} \beta_n^{n } =C(M-M')^k ,
\end{split}\end{equation*}for $k=0,1,\ldots,m$. On the other hand, for a fixed $0 <r< 1$, the left hand side above can be estimated from below as
\begin{equation*}\begin{split}&
   \int_0^1dt_1   \cdots\int_0^{t_{k-1}}dt_k\int_{|q|^2\leq t_k}\Delta u_1\wedge\dots\wedge \Delta u_k\wedge\beta_n^{n-k} \\ \geq &\int_r^1dt_1  \cdots\int_r^{t_{k-1}}dt_k\int_{|q|^2\leq r}\Delta u_1\wedge\dots\wedge \Delta u_k\wedge\beta_n^{n-k} =\frac {(1-r)^k}{k!} \int_{|q|^2\leq r}\Delta u_1\wedge\dots\wedge \Delta u_k\wedge\beta_n^{n-k} .
 \end{split}\end{equation*}
 So we get
 \begin{equation*}\begin{split}
 \int_{|q|^2\leq r}\Delta u_1\wedge\dots\wedge \Delta u_k\wedge\beta_n^{n-k} \leq \frac {Ck!(M-M')^k}{(1-r)^k} ,
 \end{split}\end{equation*}
 which implies   the local {\it Chern-Levine-Nirenberg estimate} for $QSH_m \cap C^2 $ functions.

 \begin{cor}\label{cor:uniformly-bounded} In the function class $L_M=\{u\in QSH_m(\Omega)\cap C^2(\Omega) : |u|\leq   M\},$ the integrals $\int_K\Delta u_1\wedge\dots\wedge\Delta u_k\wedge\beta_n^{n-k} $ are
uniformly bounded for any compact subset $K$, $k=1,\ldots, m$.
\end{cor}

  \begin{thm}  \label{thm:recurrence-def}
For $u_1,\dots,u_m\in QSH_m(\Omega)\cap C(\Omega)$, the recurrence relation
  \begin{equation}\label{eq:recurrence-def}
\Delta u_1\wedge\dots\wedge \Delta u_k \wedge \beta_n^{n-m}(\omega)= \int u_k\Delta u_1\wedge\dots\wedge\Delta u_{k-1}\wedge \beta_n^{n-m}\wedge \Delta \omega, \qquad k=1,\dots,m ,
  \end{equation}
for $ \omega \in  \mathcal D^{2m-2k}(\Omega)$,
   defines a  closed positive current.

   Moreover, the following weak $*$ convergence of currents of measure type holds for the standard approximations $u^t_{j}\downarrow u_j$ ($j=1,2,\dots,k$) as $t\rightarrow \infty$,
     \begin{equation}\label{eq:recurrence-convergence}
  \begin{aligned}
\Delta u_1^{t}\wedge\dots\wedge \Delta u_k^{t} \wedge \beta_n^{n-m}\rightarrow \Delta u_1\wedge\dots\wedge \Delta u_k \wedge \beta_n^{n-m}.
\end{aligned}
  \end{equation}
\end{thm}
\begin{proof} The closedness follows from definition. For $k=1$, the left hand side of \eqref{eq:recurrence-def} is the Laplace operator. The result holds.

Suppose that the result holds for $k-1$. Then $\Delta u_1\wedge\dots\wedge\Delta u_{k-1}\wedge \beta_n^{n-m} $ is a closed positive current of measure type. Thus the right hand side of \eqref{eq:recurrence-def} is well defined, and defines a linear continuous
functional on  $   \mathcal D^{2m-2k}(\Omega) $.

To show the  positivity of this current, note that  the standard approximations
  $u_j^{t} $ locally uniformly converges to
$  u_j$. Thus,
for a strongly  positive form $ \omega \in   \mathcal{D}^{2m-2k}(\Omega) $, by  the convergence \eqref{eq:recurrence-convergence} of currents of measure type for  $k-1$, we have
     \begin{equation*}
  \begin{aligned}
\Delta u_1\wedge\dots\wedge \Delta u_k \wedge \beta_n^{n-m}(\omega)&=\int u_k\Delta u_1\wedge\dots\wedge\Delta u_{k-1}\wedge\beta_n^{n-m}\wedge\Delta \omega\\
&=\lim\limits_{t\rightarrow \infty}\int u_k\Delta u^{t}_1\wedge\dots\wedge\Delta u^{t}_{k-1}\wedge\beta_n^{n-m}\wedge\Delta \omega\\
&=\lim\limits_{s\rightarrow\infty}\lim\limits_{t\to\infty}\int u^s_k\Delta u^{t}_1\wedge\dots\wedge\Delta u^{t}_{k-1}\wedge\beta_n^{n-m}\wedge\Delta \omega,
\end{aligned}
  \end{equation*}
which is nonnegative since
     \begin{equation*}\begin{split}
  \int  u^s_k\Delta u^{t}_1\wedge\dots\wedge\Delta u^{t}_{k-1}\wedge\beta_n^{n-m}\wedge\Delta \omega =\int \Delta  u^s_k\wedge \Delta u^{t}_1\wedge\dots\wedge\Delta u^{t}_{k-1}\wedge\beta_n^{n-m}\wedge \omega\geq 0,
   \end{split}  \end{equation*}by applying Stokes-type formula \eqref{eq:stokes} twice.
 Now write
  $u_k^t(q)=u_k(q)+\varepsilon_k^t(q)$. Then,
     \begin{equation*}
  \begin{aligned}
 \int & \Delta u^{t}_1\wedge\dots\wedge\Delta u^{t}_{k }\wedge\beta_n^{n-m}\wedge  \omega= \int u_k^t\Delta u^{t}_1\wedge\dots\wedge\Delta u^{t}_{k-1}\wedge\beta_n^{n-m}\wedge\Delta \omega\\
 &= \int u_k\Delta u^{t}_1\wedge\dots\wedge\Delta u^{t}_{k-1}\wedge\beta_n^{n-m}\wedge\Delta \omega+ \int \varepsilon_k^t(q)\Delta u^{t}_1\wedge\dots\wedge\Delta u^{t}_{k-1}\wedge\beta_n^{n-m}\wedge\Delta \omega  \\&
 \rightarrow \int u_k\Delta u _1\wedge\dots\wedge\Delta u _{k-1}\wedge\beta_n^{n-m}\wedge\Delta \omega=   \Delta u _1\wedge\dots\wedge\Delta u _{k }\wedge\beta_n^{n-m}(\omega),
\end{aligned}
  \end{equation*}
  by  the inductive hypothesis \eqref{eq:recurrence-convergence}  for  $k-1$ for the limit and $\varepsilon_k^t\rightarrow0$ uniformly on supp $\omega$.  Thus $\int  \Delta u^{t}_1\wedge\dots\wedge\Delta u^{t}_{k }\wedge\beta_n^{n-m}\wedge  \omega\rightarrow      \Delta u _1\wedge\dots\wedge\Delta u _{k }\wedge\beta_n^{n-m}( \omega)$ for any $\omega\in  \mathcal{D}^{2m-2k}(\Omega) $.
  By Proposition \ref{prop:positive form} and locally uniform  boundedness   of   vector measures   in Corollary \ref{cor:uniformly-bounded}, we get
  \begin{equation*}
    \left| \int_K\Delta u^{t}_1\wedge\dots\wedge\Delta u^{t}_{k }\wedge\beta_n^{n-m}\wedge  \omega\right|\leq   C_1\|\omega\|_{C(\Omega)}\int_K\Delta u^{t}_1\wedge\dots\wedge\Delta u^{t}_{k }\wedge\beta_n^{n-k} \leq CC_1\|\omega\|_{C(\Omega)}
  \end{equation*}
 where $K\supset $ supp $\omega$,   $C_1,C>0$ are absolute constants depending on $K$. We get the convergence for $(2m-2k)$-forms $\omega$ with continuous
 coefficients. Thus, \eqref{eq:recurrence-def} defines a   current of measure type.
  \end{proof}
  The measure $\Delta u^1\wedge\dots\wedge \Delta u^k \wedge \beta_n^{n-m}$ in Theorem \ref{thm:recurrence-def} is called the {\it quaternionic $m$-Hessian measure}.

Now the estimate in Theorem \ref{thm:coarea}  follows from   Lemma \ref{thm:coarea} by using Theorem \ref{thm:recurrence-def}, and
the following
  proposition also follows from   Corollary  \ref{cor:uniformly-bounded} by using Theorem \ref{thm:recurrence-def}.

\begin{prop}\label{prop:uniformly-bounded}
In the function class $L_M=\{u\in QSH_m(\Omega)\cap C (\Omega) : |u|\leq   M\},$ the  families of closed positive currents $  \Delta u_1\wedge\dots\wedge\Delta u_m\wedge\beta_n^{n-m} $ of measure type are
locally uniformly  bounded.
\end{prop}

 \begin{prop}\label{prop:sum}
  If $u,v\in C(\Omega)\cap QSH_m(\Omega)$, then $(\Delta (u+v))^m\wedge\beta_n^{n-m}\geq  (\Delta u)^m\wedge\beta_n^{n-m}+(\Delta v)^m\wedge\beta_n^{n-m}.$
  \end{prop}
\begin{proof}
Note   that if $u,v\in C^2(\Omega)\cap QSH_m(\Omega)$, we have  $(\Delta u)^i\wedge(\Delta v)^{m-i}\wedge\beta_n^m $ is positive by Proposition    \ref{prop:QSH-positivity}. So
  \begin{equation}
  \begin{aligned}
(\Delta (u+v))^m\wedge\beta_n^{n-m}&=(\Delta u)^m\wedge\beta_n^{n-m}+(\Delta v)^m\wedge\beta_n^{n-m}+\sum\limits_{p=1}^{m-1}\binom{m}{p}(\Delta u)^p\wedge(\Delta v)^{m-p}\wedge\beta_n^{n-m}\\&\geq  (\Delta u)^m\wedge\beta_n^{n-m}+(\Delta v)^m\wedge\beta_n^{n-m}.\\
\end{aligned}
  \end{equation}

  If $u,v$ is only continuous, apply the above inequality to their standard approximation $u_\epsilon , v_\epsilon $. Since $u_\epsilon, v_\epsilon$ are smooth, and $u_\epsilon\downarrow u, u_\epsilon\downarrow u, u_\epsilon+v_\epsilon\downarrow u+v$  locally uniformly. So by Theorem \ref{thm:recurrence-def}, we obtain the result  by letting $\epsilon\to 0.$
\end{proof}
It similar to  Proposition\ref{prop:uniformly-bounded} to establish the following proposition. We omit details.
\begin{prop}\label{prop:weak-compact}
In the function class $L_M=\{u\in QSH_m(\Omega)\cap C(\Omega) : |u|\leqslant M\},$  the  families of closed positive currents  $
 d_0 u_1\wedge d_1 u_1\wedge\Delta u_2\wedge\dots\wedge\Delta u_k\wedge\beta_n^{n-m}
$  of measure type are
locally uniformly  bounded.
\end{prop}
\subsection{The  comparison principle }

\begin{thm}\label{thm:compare} Let $\Omega$ be a bounded   domain  and let  $u,v\in QSH_m(\Omega)\cap C(\Omega)$. If $\{u<v\} \Subset  \Omega$, then   we have
  \begin{equation}\label{eq:compare0}
     \int_{\{u<v\}}(\Delta u  )^m\wedge\beta_n^{n-m}  \geq \int_{\{u<v\}}(\Delta v  )^m\wedge\beta_n^{n-m}
  \end{equation}

\end{thm}

 We need the following proposition to prove this theorem.
\begin{prop}\label{prop:compare}
 Let $\Omega$ be a bounded   domain with smooth boundary, and let $u ,v\in C^2(\overline{\Omega})\cap QSH_m(\Omega)$. If $u=v$ on $\partial\Omega$ and $u\leq  v$ in $\Omega$, then
 \begin{equation}\label{eq:compare}
    \int_\Omega (\Delta u  )^m\wedge\beta_n^{n-m}\geq  \int_\Omega (\Delta v  )^m\wedge\beta_n^{n-m}.
 \end{equation}
\end{prop}
\begin{proof} We can choose a   defining function $\varrho$ of $\Omega$ with  $|{\rm grad} \varrho|=1$. Then
\begin{equation}\label{eq:boundary-term}\begin{split}\int_\Omega (\Delta u)^m\wedge\beta_n^{n-m}
 &-  \int_\Omega (\Delta v)^m\wedge\beta_n^{n-m} =\int_\Omega\sum_{p=1}^{m}(\Delta v)^{p-1}\wedge\Delta(u- v)\wedge(\Delta u)^{n-p}\wedge\beta_n^{n-m}
  \\=&  \sum_{p=1}^m \int_\Omega d_0\left[d_1 \left( u- v \right)\wedge  ( \triangle
v)^{p-1} \wedge (\triangle u)^{m-p }   \wedge\beta_n^{n-m}\right ] \\=&\sum_{p=1}^m\sum_{A=0}^{2n-1} \int_{\partial\Omega} \left[d_1 \left( u- v \right)\wedge  ( \triangle
v)^{p-1} \wedge (\triangle u)^{m-p }   \wedge\beta_n^{n-m}\right ]_A
\cdot \nabla_{A0 }\varrho\, dS
\end{split}\end{equation}by using
Stokes-type formula  (\ref{lem:Stokes}). Note that we have
\begin{equation}\label{eq:boundary-term2}\begin{split}
 & \sum_{A=0}^{2n-1}  \left[d_1 \left( u- v \right)\wedge  ( \triangle
v)^{p-1} \wedge (\triangle u)^{m-p }   \wedge\beta_n^{n-m}\right ]_A
\cdot \nabla_{A0'}\varrho (q)\,\Omega_{2n}\\ =&d_0\varrho(q)\wedge d_1 \left( u-  v\right)\wedge  ( \triangle
v )^{p-1} \wedge (\triangle u )^{m-p }\wedge\beta_n^{n-m},
\end{split}\end{equation}as in \eqref{eq:boundary-term}. Since $u=v$ on $\partial\Omega$ and $u\leq v$ in $\Omega$, then for a point $q\in \partial\Omega$ with ${\rm grad} (u -v)(q)\neq 0$,   we can write  $u-v= h \varrho$ in a neighborhood of $q$ for some positive smooth function $h$. Consequently, we have ${\rm grad} (u -v)(q)= h(q) {\rm grad} \varrho$, and so  $\nabla_{A1 }(u -v)(q)= h(q)  \nabla_{A1 }\varrho (q)$ on $\partial\Omega$. Thus,
\begin{equation*}
   d_0\varrho(q)\wedge d_1 \left(u -v  \right)(q)=h(q)  d_0\varrho (q)\wedge d_1 \varrho(q) \qquad \text {on the boundary},
\end{equation*}
which is elementary   strongly positive  by Proposition \ref{prop:wedge-positive}. Since $( \triangle
v )^{p-1} \wedge (\triangle u )^{m-p }\wedge\beta_n^{n-m}$ is   also   positive  by Proposition \ref{prop:QSH-positivity}, we find that
 the right hand of  (\ref{eq:boundary-term2}) is a positive $2n$-form by definition. So the integrant in the right hand of   (\ref{eq:boundary-term}) on $\partial\Omega$  is nonnegative if ${\rm grad} (v-u)(q)\neq 0$. While if ${\rm grad} (v-u)(q)= 0$,   the integrant at $q $  in (\ref{eq:boundary-term}) vanishes. Therefore the difference in (\ref{eq:boundary-term}) is  nonnegative.
\end{proof}

{\it Proof of Theorem \ref{thm:compare}}. At first, we assume that $u,v\in QSH_m(\Omega)\cap C^2(\Omega)$. Let $G_\eta
 := \{u < v -\eta\}$. Then $G:=\{u < v  \} = \cup_{\eta>0}G _\eta $  and by Sard's theorem, $G_\eta$
 are open sets with smooth boundaries
for almost all $\eta>0$. For such $\eta$, we have
\begin{equation*}
   \int_{ G_\eta }(\Delta u  )^m\wedge\beta_n^{n-m}  \geq \int_{G_\eta}(\Delta v  )^m\wedge\beta_n^{n-m}
\end{equation*}
by Proposition  \ref{prop:compare}. \eqref{eq:compare0} follows by taking limit $\eta\rightarrow 0$.

Now if $u,v\in QSH_m(\Omega)\cap C (\Omega)$, consider  the standard approximations $u _{j}\downarrow u $, $v _{j}\downarrow v $ by smooth $  QSH_m $ functions. Denote $G_p:=\{ q\in G; u < v -1/p\}$ and $G_{j,k,p}:=\{q\in G; u_{j} < v_{k} -1/p\}$.

For any open set $G'\Subset  G$ we can choose positive integers $p_0$ and $p_1$ such
that $G' \Subset  G _{p_0}
 \Subset  G_{p_1}
 \Subset  G$. Since  $u _{j}, v _{j} $
 converge locally uniformly in $G$,
there exist   $k_0$
 such that $G'\subset G_{j,k,{p_0}}\subset  G_{p_1}\Subset  G$
 for all $j,k >k_0$. Then
 \begin{equation*}
   \int_{ G_{j,k,{p_0}}}(\Delta u_j  )^m\wedge\beta_n^{n-m}  \geq \int_{G_{j,k,{p_0}}}(\Delta v_k  )^m\wedge\beta_n^{n-m}
\end{equation*}
 for all $j,k >k_0$. Consequently,
  \begin{equation*}
   \int_{  {G}_{p_1}}(\Delta u_j  )^m\wedge\beta_n^{n-m}  \geq\int_{ {G}'}(\Delta v_k  )^m\wedge\beta_n^{n-m}.
\end{equation*}
By  convergence of currents of measure type, we get
 \begin{equation*}
   \int_{ {G }}(\Delta u   )^k\wedge\beta_n^{n-m}  \geq \int_{ {G}_{p_1}}(\Delta u   )^k\wedge\beta_n^{n-m} \geq \int_{G' }(\Delta v   )^k\wedge\beta_n^{n-m}.
\end{equation*}The result follows since the  $G'\Subset  G$ is arbitrarily chosen.
\qed

\begin{prop}Let $\Omega$ be a bounded   domain with smooth boundary, and let $u ,v\in C (\overline{\Omega})\cap QSH_m(\Omega)$. Suppose that
$(\Delta u  )^m\wedge\beta_n^{n-m} \leq (\Delta v  )^m\wedge\beta_n^{n-m}$ on $\Omega$, and
  $\underline{\lim}_{q\in \Omega}(u(q)- v(q)) \geq 0$. Then $u \geq v $ in $\Omega$.
\end{prop}
\begin{proof} Assume that $v(q_0)-u(q_0) =\eta > 0$ at
some point $q_0\in \Omega$. Thus the open set $G:=\{ D:   u(q )< v(q) - \eta /4 \}$
is not empty. Then
\begin{equation*}
   G_1:=\{ D:   u(q )< v(q) - \eta /2 +\varepsilon|q-q_0|^2\}\Subset G,
\end{equation*}
    and contains $q_0$ for sufficiently small $\varepsilon>0$. By applying the  comparison principle  in Theorem \ref{thm:compare} and  Proposition \ref{prop:sum}, we get
\begin{equation*}
   \int_{ G_1}(\Delta u   )^m\wedge\beta_n^{n-m}\geq  \int_{ G_1}(\Delta   v  +\varepsilon\triangle|q-q_0|^2  )^m\wedge\beta_n^{n-m}\geq \int_{ G_1}(\Delta   v    )^m\wedge\beta_n^{n-m}+ (8\varepsilon)^m\int_{ G_1} \beta_n^{n }
\end{equation*}
which contradicts to the assumption $(\Delta u  )^m\wedge\beta_n^{n-m} \leq (\Delta v  )^m\wedge\beta_n^{n-m}$.
\end{proof}

We also need the following  proposition for several functions.
\begin{cor} \label{cor:compare}
Let $\Omega$ be a bounded   domain and let $u_j ,v_j\in C ( {\Omega})\cap QSH_m(\Omega)$. If $u_j=v_j$ outside a compact subset of  $ \Omega$, then
 \begin{equation}\label{eq:compare=}
    \int_\Omega  \Delta u_1\wedge\dots\wedge \Delta u_m\wedge\beta_n^{n-m}=  \int_\Omega \Delta v_1\wedge\dots\wedge \Delta v_m \wedge\beta_n^{n-m}.
 \end{equation}
\end{cor} \begin{proof}  If the domain has smooth boundary and $u_j ,v_j\in C^2(\overline{\Omega})\cap QSH_m(\Omega)$, this identity
is obtained as in \eqref{eq:boundary-term} by applying
\begin{equation*}
   \Delta u_1\wedge\dots\wedge \Delta u_m- \Delta v_1\wedge\dots\wedge \Delta v_m  = \sum_{p=1}^m \Delta v_1\wedge\cdots\wedge \Delta v_{p-1} \wedge\Delta (u_p- v_p)\wedge \Delta u_{p+1}\wedge  \cdots,
\end{equation*}
since there is no boundary term in this case. The general case easily follows from approximation.
\end{proof}

 \section{Quaternionic relative $m$-extremal function  and quaternionic $m$-capacity}
For a domain $\Omega$ in $\mathbb{H}^n$ and $E\subset \Omega$, let
 \begin{equation}\label{eq:U}
   \mathcal{U}( E ,\Omega): =\{  u\in  QSH_m(\Omega) , u\vert_{\Omega}\leq    0,u\vert_{E}\leq   -1 \},
 \end{equation}
and
 \begin{equation*}
    \omega (q,E,\Omega):=\sup \{u(q); u\in \mathcal{U}( E ,\Omega) \},
 \end{equation*}
whose   upper semicontinuous regularization $ \omega^* (q,E,\Omega)$  is called a {\it  relative $m$-extremal function}  of the set $E$ in $\Omega$. The {\it $\mathcal{P}_m$-capacity} is defined as
\begin{equation*}
   \mathcal{P}_m( E,\Omega):=-\int_\Omega\omega^* (q,E,\Omega) \beta_n ^{n }.
\end{equation*}

The relative extremal function has the following simple properties:

(1) {\it (monotonicity)  if $E_1\subseteq E_2$, then  $ \omega^* (q,E_1,\Omega)  \geq \omega^* (q,E_2,\Omega) $;  if $E \subseteq D_1\subset  D_2$, then  $ \omega^* (q,E ,D_1)  \geq \omega^* (q,E ,D_2) $ for $q\in D_1$}.

(2)  {\it $ \omega^* (q,E,\Omega)\equiv 0$  if and only if $E$ is $m$-polar in $\Omega$}. The proof is the same as  the complex case \cite{klimek}.

(3) {\it  Let $    \Omega =\{\varrho<0\}$ be    $m$-hyperconvex.
If   $ E\Subset  \Omega $, then
$ \omega^* (q,E,\Omega)\rightarrow 0$ as $q\rightarrow \partial\Omega$}.

Note that $M\varrho \in\mathcal{U}( E ,\Omega)$ for a suitable $M>0$ since $E\Subset \Omega$. Then $0\geq \omega^* (q,E,\Omega)\geq M\varrho $ on $\Omega$. We must have $ \omega^* (q,E,\Omega)\rightarrow 0$ as $q\rightarrow \partial\Omega$

(4)  {\it Let $    \Omega =\{\varrho<0\}$ be a strongly $m$-hyperconvex.
If   $ E\Subset  \Omega $,
then the relative $m$-extremal function $ \omega^* (q,E,\Omega)$ admits a  quaternionic $m$-subharmonic extension to a neighborhood of the
closure $\Omega$}.

By $  \omega^* (q,E,\Omega)\geq M\varrho $ on $\Omega$ as above, the quaternionic $m$-subharmonic function
\begin{equation*}
  w(q)=\left\{ \begin{array}{ll}
  \omega^* (q,E,\Omega) ,\qquad &q\in \Omega,
  \\ M\varrho,\qquad &q\notin \Omega,
   \end{array}\right.
\end{equation*}
gives an extension to a neighborhood of $\overline{\Omega}$.

A point $q_0\in   K$ is called an {\it $m$-regular point} of the compact  set $K \Subset\Omega$ if  $\omega^* (q_0,K,\Omega) = -1$. A compact set $K\Subset  \Omega$ is called {\it  $m$-regular} in $\Omega$  if each point   of $   K$ is  $m$-regular. A function $u \in  QSH_m(\Omega)$ is called {\it maximal}
if it satisfies the  {\it maximum principle} in the class $    QSH_m(\Omega)$, i.e. for any $D\Subset\Omega$, if $v\in  QSH_m(D)$ and
$\underline{\lim}_{q\in \partial D}(u(q)- v(q)) \geq 0$, then $u \geq v $ in $D$.

Since a quaternionic   $m$-subharmonic function is subharmonic by Proposition \ref{prop:QSH-m} (2),
  a regular compact set  of the classical potential theory is  $m$-regular. In general,   an $m$-regular compact
set is always    $m'$-regular if $m'>m$. Therefore, for any    compact  subset $ K$ of an open set  $U$,
there exists an  $m$-regular compact set $ E$ such that $K \subset E \Subset U$.

\begin{prop} \label{prop:maximal}  Let $K$ be an $m$-regular compact subset of of an   $m$-hyperconvex  domain $\Omega$.
  Then, (1)  relative $m$-extremal function $\omega^*(q, K, \Omega)   $ is maximal in $ \Omega\setminus K$; (2) $\omega^*( \cdot, K, \Omega)\in C(\Omega)$; (3)
  \begin{equation}\label{eq:omega=0}
    (\Delta \omega^*(q,  K, \Omega)  )^m\wedge\beta_n^{n-m}=0\qquad \text { on }\quad \Omega\setminus K.
\end{equation}
\end{prop}
\begin{proof} (1) Suppose that $\omega^*(\cdot, K, \Omega)  $ is not maximal. Then
  there exists a domain $G \Subset  D\setminus K$ and a
function $v\in   QSH_m(G)$ such that $\underline{\lim}_{q\in \partial G}(u(q)- v(q)) \geq 0$, but $v(q_0) >
\omega^*(q_0, K, \Omega)$
at some point $q_0\in G$. Since $\omega^*(q, K, \Omega)|_K\equiv-1 $, the function
\begin{equation*}
   w(q)=\left\{\begin{array}{ll}
   \max ( v(q ),\omega^*(q, K, \Omega)),\qquad &\mbox{if } q\in G,\\\omega^*(q, K, \Omega)),\qquad &\mbox{if } q\notin G,
\end{array}
\right.
\end{equation*}belongs to
   $w\in \mathcal{U}( K ,\Omega)$ by definition \eqref{eq:U}, and so $w\leq \omega^*(\cdot, K, \Omega)$. This contradicts to $w(q_0)=v(q_0) >
\omega^*(q_0, K, \Omega)$.

(2) Consider $\Omega_j:=\{q\in\Omega; \omega^*(q , K, \Omega)<-1/j\}$ for positive integers $j$. Then $\Omega_j\subset \Omega_{j+1}$ and $\Omega_j\Subset\Omega$ since  $\Omega$ is $m$-hyperconvex. Fixed a $j_0$,  the relative $m$-extremal function can be approximated on  $\overline{\Omega}_{j_0}$ by smooth $QSH_m   $
  functions $v_{t}\downarrow  \omega^* (\cdot ,E,\Omega)$. Applying Hartogs' Lemma for subharmonic functions twice to this sequence, we  see that there exists $t_0$ such that for $t> t_0$, we have $v_t\leq0$ on  $\overline{\Omega}_{j_0}$ and simultaneously, $v_t\leq -1+ 1/j_0$ on  $K$. Then the function
  \begin{equation*}
   \widetilde{w}(q)=\left\{\begin{array}{ll}
   \max ( v_t(q )-1/j_0,\omega^*(q, K, \Omega)),\qquad &\mbox{if } q\in {\Omega}_{j_0},\\ \omega^*(q, K, \Omega)),\qquad &\mbox{if } q\notin {\Omega}_{j_0},
\end{array}
\right.
\end{equation*} belongs to
$\mathcal{U}( K ,\Omega)$, and so
\begin{equation*}
   \omega^*(q, K, \Omega) -1/j_0\leq v_t(q )-1/j_0\leq  \widetilde{w}(q)\leq \omega^*(q, K, \Omega)
\end{equation*}
 for $q\in \overline{\Omega}_{j_0}$. Consequently, $v_t $ converges uniformly to $  \omega^*(\cdot, K, \Omega)$ on compact
subsets of $\Omega$. So it is continuous.

 (3) Suppose $ (\Delta \omega^*(q,  K, \Omega)  )^m\wedge\beta_n^{n-m}$ does not vanish on $\Omega\setminus K$. There exists a ball $B(q_0,r)$ where $  (\Delta \omega^*(q,  K, \Omega)  )^m\wedge\beta_n^{n-m} \not\equiv 0$.
Let $v(q)
$ be the  Bremermann-Perron solution to the generalized Dirichlet
problem
$(\Delta v)^{m}\wedge  \beta_n ^{n-m} = 0$ on the ball with continuous boundary value $ \omega^*(\cdot,  K, \Omega)|_{\partial B(q_0,r)} $. Such a solution exits, and is unique and  continuous. The proof is exactly as  in the complex case \cite{blocki}. We omit details. It is
is maximal by construction, i.e. $v \geq \omega^*( \cdot,  K, \Omega)$ on $B(q_0,r)$.   But $v\not\equiv\omega^*(q,  K, \Omega)$,  since $  (\Delta \omega^*(q,  K, \Omega)  )^m\wedge\beta_n^{n-m} \not\equiv 0$ on $B(q_0,r)$.
 Therefore, $v(q') > u(q')$ for some $ q' \in B(q_0,r)$. But
 \begin{equation*}
    w(q)=\left\{ \begin{array}{ll}
  \omega^* (q,E,\Omega) ,\qquad &q\in \Omega\setminus B(q_0,r),
  \\ \max\{v(q ),\omega^* (q,E,\Omega)\},\qquad &q\in B(q_0,r),
   \end{array}\right.
 \end{equation*}belongs to
 $ \mathcal{U}( K ,\Omega) $. Then $w(q') > u(q')$ contradicts to the
  maximality of $\omega^*(q,  K, \Omega)$ in (1).
\end{proof}
  \subsection{Quaternionic $m$-capacity} See \cite[Section 3]{sadullaev1} for complex $m$-capacity.
Given a compact set $K$ in a domain $\Omega\subset \mathbb{H}^n$, let
\begin{equation}
\begin{aligned}
\mathcal{U}^*(K,\Omega)=\left\{u\in QSH_m(\Omega)\cap C(\Omega), u\vert_{K}\leq   -1 , \varliminf_{q\to \partial\Omega} u(q)\geq  0    \right \}.
\end{aligned}
\end{equation}
The   {\it quaternionic $m$-capacity} of the condenser $(K,\Omega)$ is defined as
\begin{equation}\label{capacity defi-K}
\begin{aligned}
C_m(K)=\inf\left\{\int_{\Omega}(\Delta u)^m\wedge\beta_n^{n-m}: u\in\mathcal{U}^*(K,\Omega)\right\}
\end{aligned}
\end{equation}
and the quaternionic $m$-capacity of an open set $U\subset\Omega $ is
\begin{equation*}
   C_m(U) = \sup\{C_m(K); K\subset U\}  .
\end{equation*}The
{\it  exterior $m$-capacity } of a set $E\subset\Omega$ is defined as
\begin{equation*}
   C_m^*(E) = \sup\{ C_m(U);  \mbox{ open   } U\supset E\}
\end{equation*}
 $m$-capacity is obviously monotonic  by definition.
\begin{prop}\label{prop:capacity}
Let $\Omega$ be a   $m$-hyperconvex domain in $\mathbb{H}^n$. Then,

(1)
For any $m$-regular compact set $K \subset \Omega$,
\begin{equation}\label{eq:capacity-K}
   C_m(K)= \int_{K}(\Delta \omega^*(q,K,\Omega))^m\wedge\beta_n^{n-m}.
\end{equation}

(2)  For any compact subset  $K \subset \Omega$, $C_m(K)=\inf\{C_m(E); \Omega \supset E\supset K$ and $ E $ is an $m$-regular   compact set  $\}$. In particular, $C_m^*(K)=C_m(K)$.

 (3) If $K$ is an $m$-regular compact subset, then
\begin{equation}\label{eq:capacity-K-2}
   C_m(K)= \sup  \left\{\int_{K} \Delta u_1  \wedge\cdots\wedge  \Delta u_m \wedge \beta_n^{n-m}; u_j\in QSH_m(\Omega)\cap C(\Omega), -1\leq u_j< 0\right\}.
\end{equation}

(4) Suppose that $\Omega$ is strongly $m$-hyperconvex.   If $U\subset \Omega$ is an open set, then
\begin{equation}\label{capacity defi}
\begin{aligned}
C_m(U)&=\sup\left\{\int_{U}(\Delta u)^m\wedge\beta_n^{n-m}:u\in QSH_m(\Omega)\cap C(\Omega), -1\leq   u<0\right\}\\
&=\sup\left\{\int_{U}(\Delta u)^m\wedge\beta_n^{n-m}:u\in QSH_m(\Omega)\cap C^{\infty}(\Omega), -1\leq   u<0\right\}
\end{aligned}
\end{equation}

(5) The exterior capacity   is monotonic, i.e.  if $E_1\subseteq E_2$, then  $C_m^*(E_1)\subseteq C_m^*(E_2)$, and
countably subadditive, i.e.  $C_m^*(\cup_j E_j)\leq \sum_j C_m^*(  E_j)$.

(6) If $U_1\subset U_2\subset\dots$ are open subsets of $\Omega$, then
$
C_m\left(\bigcup_{j=1}^\infty U_j,\Omega\right)=\lim\limits_{j\to\infty} C_m(U_j,\Omega).
$

(7) If $E\subset  D \subset \Omega$, then $C_m^*(E,D)\leq C_m^*(E,\Omega)$.
\end{prop}
\begin{proof} (1) For $u\in\mathcal{U}^*(K,\Omega)$ and any $0<\varepsilon <1$, consider the open set
 \begin{equation*}
    O:= \left\{q\in \Omega; u(q)<(1-\varepsilon)\omega^*(q,K,\Omega)-\varepsilon/2 \right\}\Subset \Omega.
 \end{equation*}
  Note that $O\supset K$. Then, we have
\begin{equation*}\begin{split}(1-\varepsilon)^m\int_{K}(\Delta \omega^*(q,K,\Omega))^m\wedge\beta^{n-m}=
& (1-\varepsilon)^m\int_{O}(\Delta \omega^*(q,K,\Omega))^m\wedge\beta_n^{n-m}\\
   \leq&\int_{O}(\Delta u)^m\wedge\beta^{n-m}\leq\int_{\Omega}(\Delta u)^m\wedge\beta_n^{n-m}
 \end{split}\end{equation*}
by the  comparison principle   and \eqref{eq:omega=0}. Letting $\varepsilon\rightarrow0$, we see that  the infimum on the right hand side of \eqref{capacity defi-K} is attained by  the relative $m$-extremal function $ \omega^*(q,K,\Omega) $.

(2) $C_m(K)\leq C_m(E)$ by monotonicity. Conversely, for   any $0<\varepsilon <1$, choose  $u\in\mathcal{U}^*(K,\Omega)$ such that $ \int_{\Omega}(\Delta u)^m\wedge\beta^{n-m}<C_m(K)+\varepsilon$. Since $\{q\in\Omega; u(q)<-1+\varepsilon\}$ is  a neighborhood of the compact set $K$, there exists an $m$-regular compact set $E$    such that $ K\subset  E\Subset  U$.  Consider
\begin{equation*}
   O:=\left\{q\in \Omega; u(q)<(1- 2\varepsilon)\omega^*(q,E,\Omega) \right\}  .
\end{equation*}Then, $ E\subset O\Subset  \{q\in\Omega; u(q)<-1+\varepsilon\}$,
 and so
\begin{equation*}\begin{split}C_m( E)&= \int_{E}(\Delta \omega^*(q,E,\Omega))^m\wedge\beta_n^{n-m}\leq \int_{O}(\Delta \omega^*(q,E,\Omega))^m\wedge\beta_n^{n-m}\\
&\leq\frac 1{ (1-2\varepsilon)^{m}}\int_{O}(\Delta u)^m\wedge\beta_n^{n-m}
    \leq \frac 1{ (1-2\varepsilon)^{m}}\int_{\Omega}(\Delta u)^m\wedge\beta_n^{n-m}\leq \frac {C_m(K)+\varepsilon}{ (1-2\varepsilon)^{m}},
 \end{split}\end{equation*}by using  \eqref{eq:capacity-K} for the $m$-regular  compact subset $E$ and the  comparison principle. The result follows by letting $\varepsilon\rightarrow0$.

(3) $ C_m(K)$ is less that or equal to  the right hand side of  \eqref{eq:capacity-K-2} by using  \eqref{eq:capacity-K}. On the other hand, for any $u_j\in QSH_m(\Omega)\cap C(\Omega)$ with $ -1\leq u_j< 0$, consider
\begin{equation*}
   v_j(q):=\max\left\{(1+\varepsilon)\omega^*(q,K,\Omega), \frac {u_j(q)-\varepsilon/2}{1+\varepsilon/2}\right\}.
\end{equation*}
  Then,
$v_j\in QSH_m(\Omega)\cap C(\Omega)$ with $ -1\leq v_j< 0$, $\lim_{q\rightarrow\partial\Omega}v_j(q)=0$, and $v_j\equiv(1+\varepsilon)\omega^*(\cdot,K,\Omega)$ near the boundary. We get
\begin{equation*}\begin{split} (1+\varepsilon)^{ m}\int_{\Omega}(\Delta \omega^* )^m\wedge\beta_n^{n-m}= \int_{\Omega} \Delta v_1  \wedge\cdots\wedge  \Delta v_m \wedge \beta_n^{n-m}\geq \frac {1}{(1+\varepsilon/2)^m}\int_{K} \Delta u_1  \wedge\cdots\wedge  \Delta u_m \wedge \beta_n^{n-m}.
 \end{split}\end{equation*}by using Corollary \ref{cor:compare} and $v_j\equiv(u_j -\varepsilon/2)/(1+\varepsilon/2)$ in a neighborhood of $K$. Letting $\varepsilon\rightarrow0$, we get the another direction of inequality,
  since $(\Delta \omega^* )^m\wedge\beta^{n-m}=0$ on $\Omega\setminus  K $.

  (4) For any $u\in QSH_m(\Omega)\cap C(\Omega)$ with $-1\leq   u<0$, we have $C_m(U) \geq C_m(K)\geq  \int_{K}(\Delta u)^m\wedge\beta_n^{n-m}$ by (3). Then $C_m(U)  \geq  \int_{U}(\Delta u)^m\wedge\beta_n^{n-m}$ since $K$ can be arbitrarily chosen. Thus $C_m(U)$ is larger than or equal to the right hand side of \eqref{capacity defi}.

Since $\Omega$ is a strongly $m$-hyperconvex domain, the relative $m$-extremal function $ \omega^* (q,E,\Omega)$ admits an quaternionic $m$-subharmonic extension to a neighborhood of the
closure $\Omega$, and so it can be approximated in a neighborhood $U$ of $\overline{\Omega}$ by $QSH_m \cap C^\infty $
  functions $v_{j}\downarrow  \omega^* (q,K,\Omega)$. Hence,
  \begin{equation*}\begin{split}C_m( K)&= \int_{K}(\Delta \omega^*(q, K ,\Omega))^m\wedge\beta_n^{n-m}= \int_{\Omega}(\Delta \omega^*(q, K ,\Omega))^m\wedge\beta_n^{n-m} \leq \overline{\lim} _{j\rightarrow \infty} \int_{\Omega}(\Delta v_{j} ))^m\wedge\beta_n^{n-m}\\
&\leq \overline{\lim}_{j\rightarrow \infty} (1+\varepsilon)^{ m}\int_{\Omega}(\Delta w_{j} ))^m\wedge\beta_n^{n-m}
 \end{split}\end{equation*}
if we denote $ w_{j}= (v_{j}-\varepsilon)/(1+\varepsilon)$. Here $-1\leq w_{j}<0$  if $j$  is large. So $C_m( K)$ is controlled by the right hand side of  \eqref{capacity defi} multiplying $(1+\varepsilon)^{ m}$. The result follows by letting $\varepsilon\rightarrow0$.

  (5) The monotonicity of $C_m^*(E ) $ follows from the monotonicity of $C_m ( K)$  for compact sets $K$. If $E_j$'s are open sets, then
  \begin{equation*}\begin{split}
  C_m( \cup_j E_j)&=
 \sup\left\{\int_{\cup_j E_j}(\Delta u)^m\wedge\beta_n^{n-m}:u\in QSH_m(\Omega)\cap C (\Omega), -1\leq   u<0\right\} \\
 & \leq\sup\left\{\sum_j\int_{ E_j}(\Delta u)^m\wedge\beta_n^{n-m}:u\in QSH_m(\Omega)\cap C (\Omega), -1\leq   u<0\right\}\leq \sum_j C_m(  E_j).
 \end{split}\end{equation*}
  In general, we find an open set $U_j\supset E_j$ such that $C_m( U_j) -C_m^*(  E_j)\leq \varepsilon/2^j$. Then
  \begin{equation*}
     \sum_j C_m^*(  E_j)\geq \sum_j C_m( U_j)-\varepsilon \geq C_m( \cup_j U_j)-\varepsilon\geq C_m( \cup_j E_j)-\varepsilon.
  \end{equation*}
  We get the result by letting $\varepsilon\rightarrow0$.

   (6) It is obvious by definition.
 \end{proof}

By (4) and (5), we get a useful  estimate: for  a strongly $m$-hyperconvex domain $\Omega$, there exists a neighborhood $\Omega'\supset\overline{\Omega}$ such that
\begin{equation}\label{eq:outside}
\int_{U}(\Delta u_1) \wedge\cdots\wedge (\Delta u_m)\wedge \beta^{n-m}  \leq C_m(U)
\end{equation} for any $u_j\in QSH_m(\Omega)\cap C(\Omega)$ with $ -1\leq u_j< 0$ on $\Omega$ and $|u_j|\leq1$ on $\Omega'$.
\begin{prop}     If $E \subset B(0,r)  $, $r < 1$, then
\begin{equation}\label{eq:C-P}
   C_m^*(E,B)\leq \frac {m!\mathcal{P}_m(E,B)}{(1-r^2)^m}.
\end{equation}
\end{prop}
 \begin{proof}   It is sufficient to prove \eqref{eq:C-P} for $m$-regular compact set $E$.  Apply Theorem   \ref{thm:coarea} for $\varrho(q)=|q|^2-1$, $\Omega=B$ and $u=\omega= \omega (q,E,B)$ repeatedly to get
 \begin{equation*}
    \int_0^1dt_1  \cdots\int_0^{t_{m-1}}dt_m\int_{|q|^2\leq t_m}(\Delta \omega)^m\wedge\beta^{n-m}\leq   \int_0^1dt_1  \int_{|q|^2\leq t_1} \Delta \omega \wedge\beta^{n-1}= - \int_{B}\omega  \beta_n^{n }=\mathcal{P}_m( E,\Omega).
 \end{equation*}
 On the other hand,
 \begin{equation*}\begin{split}
   \int_0^1dt_1   \cdots\int_0^{t_{m-1}}dt_m\int_{|q|^2\leq t_m}(\Delta \omega)^m\wedge\beta_n^{n-m}&\geq \int_{r^2}^1dt_1  \cdots\int_{r^2}^{t_{m-1}}dt_m\int_{|q|^2\leq{r^2}}(\Delta \omega)^m\wedge\beta^{n-m}\\&=\frac {(1-r^2)^m}{m!} \int_{|q|^2\leq r^2}(\Delta \omega)^m\wedge\beta_n^{n-m}.
 \end{split}\end{equation*}The estimate follows.
 \end{proof}

 \section{ The quasicontinuity of  quaternionic $m$-subharmonic functions and the Bedford-Taylor theory}

\begin{lem}\label{lem: inner-product} \cite[Corollary 3.1]{WZ}
If $u,v\in C^2(\Omega)$ and let $\alpha$ be a positive (2n-2)-form. Then
\begin{equation}
\begin{aligned}
\left|\int_{\Omega}d_0 u\wedge d_1 v\wedge\alpha\right|^2\leq  \int_{\Omega}d_0 u\wedge d_1 u\wedge\alpha\cdot \int_{\Omega}d_0 v\wedge d_1 v\wedge\alpha.
\end{aligned}
\end{equation}
\end{lem}

  \begin{thm}\label{thm:quasicontinuity}
Any bounded  quaternionic $m$-subharmonic function is continuous almost everywhere with respect to $m$-capacity, i.e., given $u\in QSH_m(\Omega)$ and any $\epsilon>0$, there exists an open set $U\subset \Omega$ such that $C_m(U,\omega)<\epsilon$ and $u$ is continuous on $\Omega\setminus U$.
\end{thm}
\begin{proof}
Firstly, we establish an integral inequality for $QSH_m$ functions on $B .$
Let $\mathcal{L}$ be the class of smooth $QSH_m$ functions $u$ on the ball $B(0, 1 + \delta)$ for  $\delta> 0$, such
that $|u| \leq 1$.
Consider functions $v,u, u_1,\dots,u_m\in \mathcal{L} $ such that $\varphi_0=v-u\geq  0$ in $B$ and $\varphi_0=const$ on the sphere $S=\partial B$. Then if we denote $\Theta:=d_1 u_1\wedge\Delta u_2\wedge\dots\wedge\Delta u_m\wedge\beta_n^{n-m}$,  we get
\begin{equation}
\begin{aligned}
\int_B \varphi_0\Delta u_1\wedge\dots\wedge\Delta u_m\wedge\beta_n^{n-m}&=\int_B \varphi_0d_0\Theta
 =\sum_{A=0}^{2n-1} \int_S \varphi_0 \Theta_A \tau( \mathbf{{n}})_{A0}dS-\int_ B d_0\varphi_0\wedge \Theta \\&=\varphi_0\sum_{A=0}^{2n-1} \int_S \Theta_A \tau( \mathbf{{n}})_{A0}dS-\int_ B d_0\varphi_0\wedge \Theta \\&=\varphi_0 \int_B  d_0\Theta-\int_ B d_0\varphi_0\wedge \Theta \\
&\leq   C \|\varphi\|_S-\int_ B d_0\varphi_0\wedge d_1 u_1\wedge\dots\wedge\Delta u_m\wedge\beta_n^{n-m} \\
\end{aligned}
\end{equation}by using Stokes-type formula \eqref{eq:stokes} to functions in $ \mathcal{L}$,
  where $ C $ is  an absolute constant  independent of $u_1,\dots,u_m\in \mathcal{L} $ by Corollary \ref{cor:uniformly-bounded}.
Applying Lemma  \ref{lem: inner-product}    to $u=\varphi, v=u_1$ and closed positive form $\alpha=\Delta u_2\wedge\dots\wedge\Delta u_m\wedge\beta_n^{n-m}  $, and using Stokes-type formula \eqref{eq:stokes} twice, we get
\begin{equation}\label{eq:phi0}
\begin{aligned}
\left|\int_ B d_0\varphi_0\wedge d_1 u_1\wedge\alpha\right|^2&\leq  \left(\int_B d_0 u_1 \wedge d_1 u_1\wedge\alpha\right)\left(\int_B d_0 \varphi_0 \wedge d_1 \varphi_0\wedge\alpha\right)\\
 &\leq      C\left(\varphi_0\sum_{A=0}^{2n-1} \int_S\left(d_1\varphi_0\wedge\alpha\right)_A \tau( \mathbf{{n}})_{A0}dS+\int_B \varphi_0\Delta\varphi_0\wedge\alpha\right) \\&=   C\left(\varphi_0\int_B \triangle \varphi_0\wedge\alpha+\int_B \varphi_0\Delta\varphi_0\wedge\alpha\right)  \\
&\leq   C   \left(2C \|\varphi_0\|_S+\int_B 2\left(\varphi_0\Delta\left(\frac{u+v}{2}\right)-\varphi_0\Delta v\right)\wedge\alpha\right)\\
&\leq  C\left(2C \|\varphi_0\|_S+2\int_B \varphi_0\Delta \varphi_0^{+}\wedge\alpha\right),
\end{aligned}
\end{equation}
where $\varphi_0^{+}=\frac{u+v}{2}\in \mathcal{L}$. The   second inequality follows from locally uniform estimate in Proposition \ref{prop:weak-compact}, $\varphi_0\vert_S=\|\varphi_0\|_S$ and $$ \left|\int_B\Delta \varphi_0\wedge\alpha\right|\leq   \left|\int_B\Delta(u+v)\wedge\alpha\right|\leq   2C   ,$$
while the last inequality in \eqref{eq:phi0} follows from the fact $\varphi_0\geq  0$ and $\Delta v\wedge\alpha\geq  0$.

Applying this procedure repeatedly, we obtain the inequality
\begin{equation}\label{eq:kappa}
\begin{aligned}
\int_B\varphi_0 \Delta u_1\wedge\dots\wedge\Delta u_m\wedge\beta_n^{n-m}\leq   \gamma\left(\|\varphi_0\|_S+\int_B \varphi_0(\Delta \varphi_0^{+})^m\wedge\beta_n^{n-m}\right)^\kappa,
\end{aligned}
\end{equation}
for some absolute constants $\gamma ,\kappa>0$.

Since the capacity is countably subadditive, it suffices to prove the theorem for the unit $B\subset \Omega$ and show that for any $\epsilon >0$ there exists an open set $U\subset B'$ such $C_m(U\cap B',B)<\epsilon$ and $u$ is continuous in $B'\setminus U,$ where $B'=B(0,\frac{1}{2})$. Assume $-1\leq   u\leq   0$. If replace $u$ by $\max\{u(q),v(q))\}$ with $v(q)=2(|q|^2-\frac{3}{4}),$ then $v(q)\vert_{\partial B}=\frac{1}{2}>0>u(q)$, i.e. $u\equiv v$ in a neighborhood of the sphere $S=\partial B$. Let  $u_p \downarrow u$, $v_p\downarrow  v$ be the standard approximations. Note that $u_p\equiv v_p$ in in a neighborhood of $S$ for $p>p_0$.
We can assume the sequence $\int_B u_p(\Delta u_p)^m\wedge\beta_n^{n-m}$ has a limit by passing to subsequence if necessary, since it is bounded by Proposition \ref{prop:uniformly-bounded}.  For a fixed $\sigma>0$, consider $U_{p,N}(\sigma):=\{q\in B' : u_p(q)-u_{p+N}(q)>\sigma\}$, then we have $
 U_{p,N}(\sigma)\subset U_{p,N+1}(\sigma),$ and
$\bigcup_{N=1}^\infty U_{p,N}=U_p(\sigma):=\{q\in B':u_p(q)-u(q)>\sigma\}.
$
Then we have
  \begin{equation}
  C_m\left(\bigcup_{N=1}^\infty U_{p,N}(\sigma)\right)=C_m(U_p(\sigma))=\lim\limits_{N\to\infty}C_m(U_{p,N}(\sigma))
  \end{equation}
by    Proposition \ref{prop:capacity} (6).

Denote $\varphi_{p,N}:=u_p-u_{p+N}.$ Since the open set $ U_{p,N}(\sigma)\subset B'\Subset B$, it follows from   \eqref{eq:outside}  that
  \begin{equation}
  \begin{aligned}
C_m(U_{p,N}(\sigma))&=\sup\left\{\int_{U_{p,N}(\sigma)}(\Delta u)^m\wedge\beta_n^{n-m}:u\in \mathcal{L}\right\}
\\& \leq   \sup\left\{\frac{1}{\sigma}\int_{U_{p,N}(\sigma)}\varphi_{p,N}(\Delta u)^m\wedge\beta_n^{n-m}:u\in \mathcal{L}\right\}\\
&\leq   \sup\left\{\frac{1}{\sigma}\int_{B}\varphi_{p,N}(\Delta u)^m\wedge\beta_n^{n-m}:u\in \mathcal{L}\right\}\\
&\leq  \frac{\gamma}{\sigma}\left(\|v_p-v\|_S+\int_B \varphi_{p,N}(\Delta \varphi_{p,N}^{+})^m\wedge\beta_n^{n-m}\right)^\kappa,\\
\end{aligned}
  \end{equation}by the estimate \eqref{eq:kappa},
where $\varphi_{p,N}^+:=(u_p+u_{p+N})/2.$ Note that
   \begin{equation}
  \begin{aligned}
(\Delta \varphi_{p,N}^{+})^m\wedge\beta_n^{n-m}=2^{-m}(\Delta u_p+\Delta u_{p+N})^m\wedge\beta_n^m=2^{-m}\sum\limits_{k=0}^m \binom mk(\Delta u_p)^k\wedge(\Delta u_{p+N})^{m-k}\wedge\beta_n^m.
\end{aligned}
  \end{equation}
  It is sufficient to prove $\int_{B}(u_p-u_{p+N})(\Delta u_p)^k\wedge(\Delta u_{p+N})^{m-k}\wedge\beta_n^m$ tends to $0$  uniformly  as $N\to\infty $ and then $p\to\infty $.

  For any closed $C^2$ smooth $2(n-1)$-form $\alpha$, i.e., $d_0\alpha=0,d_1\alpha=0,$ such that $\Delta u_{p+N}\wedge\alpha\geq 0$, we have
   \begin{equation*}
  \begin{aligned}
\int_B u_p\Delta u_{p+N}\wedge\alpha
&=\sum_{A=0}^{2n-1} \int_S u_p(d_1u_{p+N}\wedge\alpha)_A \tau( \mathbf{{n}})_{A0}dS+\int_Bd_1u_{p+N}\wedge d_0 u_p\wedge\alpha\\
&=\sum_{A=0}^{2n-1} \int_S u_p(d_1u_{p+N}\wedge\alpha)_A \tau( \mathbf{{n}})_{A0}dS+\sum_{A=0}^{2n-1} \int_Su_{p+N} (d_0 u_p\wedge\alpha)_A \tau( \mathbf{{n}})_{A1}dS\\
&\quad+\int_B u_{p+N}\Delta u_p\wedge\alpha\\
&\leq    A_{p,N}+\int_B u_p\Delta u_p\wedge\alpha,
\end{aligned}
  \end{equation*} by using Stokes-type formula \eqref{eq:stokes}  and $d_0d_1=-d_1d_0$,
  where
  \begin{equation}
    A_{p,N}:=\sum_{A=0}^{2n-1} \int_S \left[v_p(d_1v_{p+N}\wedge\alpha)_A \tau( \mathbf{{n}})_{A0}+v_{p+N} (d_0 v_p\wedge\alpha)_A \tau( \mathbf{{n}})_{A1}\right]dS,
  \end{equation}since $u_p=v_p$ in a neighborhood of $S$ for   $p>p_0$.
  Similarly,
     \begin{equation}
  \begin{aligned}
 \int_B u_{p+N}\Delta u_p\wedge\alpha &=B_{p,N}+\int_B u_p  \Delta u_{p+N} \wedge\alpha
 \geq  B_{p,N}+\int_B u_{p+N}\Delta u_{p+N}\wedge\alpha,
\end{aligned}
  \end{equation}
by $u_p\geq u_{p+N}$ and $ \Delta u_p\wedge\alpha\geq 0$,  where
  \begin{equation*}
     B_{p,N}:=\sum_{A=0}^{2n-1} \int_S\left [v_{p+N}(d_1v_p\wedge\alpha)_A \tau( \mathbf{{n}})_{A0}dS+ v_p (d_0 v_{p+N}\wedge\alpha)_A \tau( \mathbf{{n}})_{A1}\right]dS.
  \end{equation*}

   Repeating this procedure, finally we get
     \begin{equation}\label{eq:difference-pN}
  \begin{aligned}
 &\int_{B}(u_p-u_{p+N})(\Delta u_p)^k\wedge(\Delta u_{p+N})^{m-k}\wedge\beta_n^m  \\
\leq &  \sigma(v,p, N)+\int_B u_p(\Delta u_p)^m\wedge\beta_n^{n-m}-\int_B u_{p+N}(\Delta u_{p+N})^m\wedge\beta_n^{n-m}
\end{aligned}
  \end{equation}where $\sigma(v,p, N)$ is the sum of terms of type $A_{p,N}$ and $B_{p,N}$ above.
  Because the sequence $\{v_p\}$ converges in the $C^2(\overline{B})$, we have
    \begin{equation}
   \begin{split}
 A_{p,N}&\rightarrow \sum_{A=0}^{2n-1} \int_S \left[v_p(d_1v_{p }\wedge\alpha)_A \tau( \mathbf{{n}})_{A0}+v_{p } (d_0 v_p\wedge\alpha)_A \tau( \mathbf{{n}})_{A1}\right]dS
 \\&= \int_B (d_0 v_p \wedge d_1v_{p }+ v_p   d_0 d_1v_{p }) \wedge\alpha+\int_B ( d_1v_{p }\wedge d_0 v_p+ v_p d_1  d_0v_{p })\wedge\alpha=0,
 \end{split}
  \end{equation}
as $N\to\infty$, by using Stokes-type formula \eqref{eq:stokes} again. Similarly  $B_{p,N}\rightarrow 0$ as $N\to\infty$.
Since the sequence $\int_B u_p(\Delta u_p)^m\wedge\beta_n^{n-m}$ has a limit as  $p\rightarrow \infty$,  the right hand side of \eqref{eq:difference-pN} tends to $0$. Hence   $$\lim\limits_{p\to\infty}C_m(U_p(\sigma))=\lim\limits_{p\to\infty}\lim\limits_{N\to\infty}(U_{p,N}(\sigma))=0.$$

 Now for fixed $\epsilon>0$, there exist $p_j>0$ such that if we denote $U_{p_j}:=U_{p_j}( {1}/{j})$ for $\sigma=\frac{1}{j}$,   we have $C_m(U_{p_j})\leq   \frac{\epsilon}{2^j}$. Since $u_p(q)-u(q)<\frac{1}{j}$ for $p>p_j$ outside the set
 $U_{p_j}$, then we see that $u_p$ convergence to $u$ uniformly outside the open set $U=\cup_{j=1}^{\infty}U_{p_j}$. Since $u_p\in C^{\infty}(B),$   $u$ is continuous outside $U$, and
  $$C_m(U )=C_m\left(\bigcup_{j=1}^\infty U_{p_j}\right)\leq  \sum\limits_{j=1}^\infty C_m(U_{p_j})\leq   \epsilon .$$
 The theorem is proved. \end{proof}

\begin{prop}
Let   $u_1,\dots u_m\in QSH_m(\Omega)\cap L_{loc}^{\infty}(\Omega).$ Then,
(1) the recurrence relation
  \begin{equation}
  \begin{aligned}
\Delta u_1\wedge\dots\wedge\Delta u_k&\wedge\beta_n^{n-m}(\omega)=\int u_k\Delta u_1\wedge\dots\wedge\Delta u_{k-1}\wedge\beta_n^{n-m}\wedge\Delta \omega
\end{aligned}
  \end{equation}
  for $\omega\in   \mathcal{D}^{2m-2k}(\Omega), k=1,\dots,m,$ defines a closed positive $ 2(n-m+k) $-current.

(2) The following convergence of closed positive currents (of measure type)  holds for the standard approximations $u_i^t\downarrow u_i, i=1,\dots,m,$ as $t\to\infty:$
  \begin{equation}
  \begin{aligned}\label{LLLLL}
\Delta u_1^t\wedge\dots\wedge\Delta u_k^t&\wedge\beta_n^{n-m}\mapsto\Delta u_1\wedge\dots\wedge\Delta u_k\wedge\beta_n^{n-m}, \quad k=1\dots ,m.
\end{aligned}
  \end{equation}
\end{prop}
\begin{proof}
Let us prove the theorem by induction on $k$. The case $k=1$   is obvious.

Assume that it holds for $k-1$. Then for a fixed strongly positive form $\omega\in \mathcal{D}^{2m-2k}(\Omega)$, we have
  \begin{equation}
  \begin{aligned}
\int u_k^s\Delta u_1^t\wedge\dots\wedge\Delta u_{k-1}^t\wedge\beta_n^{n-m}\wedge\Delta \omega=
\int \Delta u_1^t\wedge\dots\wedge\Delta u_{k-1}^t\wedge\Delta u_k^s\wedge\beta_n^{n-m}\wedge\omega\geq  0,
\end{aligned}
  \end{equation}by Proposition  \ref{prop:QSH-positivity} for smooth $  QSH_m  $,
  which yields the limits $\int u^s_{k}\Delta u_1\wedge\dots\wedge\Delta u_{k-1}\wedge\beta_n^m\wedge\Delta\omega\geq  0$   as $t\to\infty$. If let  $s\to\infty$,
  we find that $\int u_k\Delta u_1\wedge\dots \wedge\Delta u_{k-1}\wedge\beta_n^{n-m}\wedge\Delta \omega\geq  0$. Hence, the current $\Delta u_k\wedge \Delta u_1\wedge\dots \wedge\Delta u_{k-1}\wedge\beta_n^{n-m}$ is positive. It is closed by definition.

  To prove (2), note that if the convergence
   \begin{equation}
\label{LOPO}
\mathcal{E}:=u_1^t \Delta u_2^t\wedge\dots\wedge\Delta u_{k }^t \wedge\beta_n^{n-m} -u_1\Delta u_2\wedge\dots\wedge\Delta u_{k }\wedge\beta_n^{n-m}\longrightarrow0, \qquad \text {as } \quad t\to\infty,
  \end{equation}
  is valid for $k $, then $(\ref{LLLLL})$ is valid for $k$, since
  \begin{equation*}
     \int \Delta u_1^t\wedge\dots\wedge\Delta u_k^t\wedge\beta_n^{n-m}\wedge\omega=\int u_1^t\Delta u_2^t\wedge\dots\wedge\Delta u_k^t\wedge\beta_n^{n-m}\wedge\Delta\omega ,
  \end{equation*}
 for $\omega\in \mathcal{D}^{2m-2k}(\Omega).$ So it suffices to prove $(\ref{LOPO})$ for $k $, provided that $(\ref{LLLLL})$ is valid for $k-1$.

 By the quasicontinuity in Theorem \ref{thm:quasicontinuity}, for a fixed $\epsilon>0$, we can find an open $U\subset\Omega$ such that $C_m(U)<\epsilon$ and $u_1\in C(\Omega\setminus U)$. Let $\widetilde{u}\in C(\Omega)$ satisfy $u_1\equiv \widetilde{u}$ on   $\Omega\setminus U$ and $\|\widetilde{u}\|_{\Omega}\leq   \|u\|_{\Omega}$.   Denote $E_{\omega}:=\operatorname{  supp} \omega $. Then,
   \begin{equation*}
  \begin{aligned}
 \left|\mathcal{E}\wedge \omega\right|
\leq  &\left|\int_{E_{\omega}} (u_1^t-u_1)\Delta u_2^t\wedge\dots\wedge\Delta u_k^t\wedge\beta_n^{n-m}\wedge\omega\right|\\
&+\left|\int_{\Omega} u_1\left(\Delta u_2^t\wedge\dots\wedge\Delta u_k^t  -\Delta u_2\wedge\dots\wedge\Delta u_k \right)\wedge\beta_n^{n-m}\wedge\omega\right|\\
 \leq   &\left|\int_{E_{\omega}\setminus U}(u_1^t-u_1)\Delta u_2^t\wedge\dots\wedge\Delta u_k^t\wedge\beta_n^{n-m}\wedge\omega\right| +\left|\int_{E_{\omega}\cap U}(u_1^t-u_1)\Delta u_2^t\wedge\dots\wedge\Delta u_k^t\wedge\beta_n^{n-m}\wedge\omega\right|\\
  &+\left|\int_{E_\omega\cap U }(u_1-\widetilde{u})\left(\Delta u_2^t\wedge\dots\wedge\Delta u_k^t -\Delta u_2\wedge\dots\wedge\Delta u_k\right)\wedge\beta_n^{n-m}\wedge\omega\right|\\
  &+\left|\int_{\Omega}\widetilde{u}\left(\Delta u_2^t\wedge\dots\wedge\Delta u_k^t  -\Delta u_2\wedge\dots\wedge\Delta u_k\right)\wedge\beta_n^{n-m}\wedge\omega\right|.
  \end{aligned}
  \end{equation*}
 The integral over the sets $E_{\omega}\setminus U$ on the right hand side tends to zero as $t\to\infty$ since $u_1^t\to u_1$ uniformly in $E_{\omega}\setminus U$, while the forth  integral over $\Omega$ tends to zero because
$$\lim\limits_{t\to+\infty}\Delta u_2^t\wedge\dots\wedge\Delta u_k^{t}\wedge \beta_n^{n-m}=\Delta u_2 \wedge\dots\wedge\Delta u_k\wedge\beta_n^{n-m},$$ as currents of measure type  by the assumption of induction, and
  $\widetilde{u}$
 continuous on $\Omega$. The second and third integrals reduces to estimating integrals of the type
$$\int_{E_{\omega\cap U}} \Delta v_2\wedge\dots\wedge\Delta v_k\wedge\beta_n^{n-m}\wedge\omega,$$
where $v_2,\dots,v_k\in QSH_m(\Omega)\cap L_{loc}^{\infty}(\Omega),$ which are small because the capacity $C_m(U)<\epsilon$ is small.

 At last, a positive current is   a current of measure type    by  Proposition 3.4  in \cite{wan-wang}.
\end{proof}

  \section{The fundamental solution  of the  quaternionic  $m$-Hessian operator and the  $m$-Lelong number }

\begin{prop} Let $\kappa_m=\frac{2n}{m}-1$. Then
the function $K_m(q): =\frac{-1}{ |q-a |^{2\kappa_m}}$ is $QSH_m$ and is the fundamental solution to the quaternionic   $m$-Hessian operator $\mathcal{H}_m$,
i.e.
  \begin{equation}
  \begin{aligned}\label{LLLL}
  \mathcal{H}_m(K_m)=C_{m,n}\delta _a
\end{aligned}
  \end{equation}
 where $C_{m,n}=\frac{8^mn!\pi^{2n} \kappa_m^m}{(2n)!m!(n-m)!}.$
   \end{prop}

\begin{proof}:
Without loss of generality, we may assume that $a=0$. Denote $K_{m,\epsilon}:=\frac{-1}{({|q|^2+\epsilon})^{\kappa_m}}$
 Then,
 \begin{equation}\label{eq:d1-K}
   d_1   K_{m,\epsilon}
  = \frac{\kappa_m d_1|q|^2}{(|q|^2+\epsilon)^{\kappa_m+1}},
  \end{equation}
 and
\begin{equation}
  \begin{aligned}
  \Delta K_{m,\epsilon}
  = d_0 \left(\frac{\kappa_m d_1|q|^2}{(|q|^2+\epsilon)^{\kappa_m+1}}\right)
= -\frac{\kappa_m(\kappa_m+1)}{(|q|^2+\epsilon)^{\kappa_m+2}}d_0|q|^2\wedge d_1|q|^2+\frac{8\kappa_m\beta_n}{(|q|^2+\epsilon)^{\kappa_m+1}}=:A+B.
\end{aligned}
  \end{equation}
 Hence,
  \begin{equation}
  \label{ppppp}
( \Delta K_{m,\epsilon})^p\wedge\beta_n^{n-p}=(pA\wedge B^{p-1}+B^p)\wedge\beta_n^{n-p},\qquad p= 1,\ldots,m,
   \end{equation}
  by $\omega\wedge\omega=0$ for any $1$-form $\omega$.
Now apply
  \begin{equation}
   d_0|q|^2\wedge d_1|q|^2=4 \sum_{l=0}^{n-1}|q_l|^2 \omega^l\wedge \omega^{n+l}+\sum_{|j-k|\not = n}a_{jk}\omega^j\wedge \omega^k
  \end{equation}
(cf. \cite[(3.12)]{wang2}) to  (\ref{ppppp}) to get
  \begin{equation*}
  \begin{aligned}
&{(\Delta K_{m,\epsilon})}^p\wedge  \beta_n ^{n-p}\\
                          =&\left[ -\frac{4p\kappa_m(\kappa_m+1)}{(|q|^2+\epsilon)^{\kappa_m+2}} \sum_{l=0}^{n-1}|q_l|^2 \omega^l\wedge \omega^{n+l}\wedge\left(\frac{8\kappa_m\beta_n}{(|q|^2+\epsilon)^{\kappa_m+1}}\right)^{p-1}
                              +\left(\frac{8\kappa_m\beta_n}{(|q|^2+\epsilon)^{\kappa_m+1}}\right)^{p}\right]\wedge \beta_n ^{n-p}\\
                           =& \frac{-4p(\kappa_m+1)\kappa_m^{p}(n-1)!8^{p-1}|q|^2}{(|q|^2+\epsilon)^{1+(\kappa_m+1)p}}\Omega_{2n}
+
\frac{8^{p}n! \kappa_m ^p}{{(|q|^2+\epsilon)^{(\kappa_m+1)p}}}\Omega_{2n} \\=& \frac{ 4\kappa_m ^{p }(n-1)!8^{p-1}|q|^2}{(|q|^2+\epsilon)^{1+(\kappa_m+1)p}}(- p (\kappa_m+1)+2n )\Omega_{2n}
+
\epsilon \frac{8^{p}n!  \kappa_m ^p}{{(|q|^2+\epsilon)^{1+(\kappa_m+1)p}}}\Omega_{2n}\geq 0
\end{aligned}
  \end{equation*} by $- p (\kappa_m+1)+2n=2n (1-p/m)\geq0$. Thus, $K_{m,\epsilon}\in QSH_m$ by definition, and so is $K_{m }\in QSH_m$ by $K_{m,\epsilon} \downarrow K_{m }$.
  In particular,
  \begin{equation*}
  \begin{aligned}
  (\Delta K_{m,\epsilon}) ^m\wedge  \beta_n ^{n-m}
= \epsilon \frac{8^{m}n!\kappa_m ^m}{{(|q|^2+\epsilon)^{2n+1}}}\Omega_{2n}.
\end{aligned}
  \end{equation*}
  Letting $\epsilon\to 0$, we get
  \begin{equation}
\left( \Delta K_m\right)^m\wedge  \beta_n ^{n-m}=0 \qquad \text{ on } \quad  \mathbb{H}^n\setminus \{0\}.
  \end{equation}

    For any $\varphi \in \mathbb{C}_0^{\infty}(\mathbb{R}^{4n})$, by rescaling   $q=q'\epsilon^{\frac{1}{2}}$,
  we get
   \begin{equation*}
  \begin{aligned}
\lim_{\epsilon\to0} \int_{\mathbb{R}^{4n}} \frac{\epsilon  }{{(|q|^2+\epsilon)^{2n+1}}}\varphi(q)dV(q)
 =\lim_{\epsilon\to0} \int_{\mathbb{R}^{4n}} \frac{\varphi(q'\epsilon^{\frac{1}{2}})}{(\|q'\|^2+1)^{2n+1}}dV(q')=\frac{S_{4n}}{4n} \varphi(0) ,
  \end{aligned}
  \end{equation*}
  by
   \begin{equation*}
  \begin{aligned} \label{iiiiiiiii}
\int_{\mathbb{R}^{4n}}\frac{1}{(|q|^2+1)^{2n+1}}dV(q)&=\lim_{R\to\infty}S_{4n}\int_{0}^{R}\frac{r^{4n-1}}{(1+r^2)^{2n+1}}dr =\lim_{R\to\infty}S_{4n}\int_{0}^{\arctan R}\frac{\tan^{4n-1}\theta}{\sec^{4n}\theta}d\theta \\
&=   \lim_{R\to\infty}S_{4n}\int_{0}^{\arctan R}\sin^{4n-1}d\sin\theta
=\lim_{R\to\infty}S_{4n}\int_{0}^{\frac{R}{\sqrt{1+R^2}}}t^{4n-1}dt\\
& =\lim_{R\to\infty}S_{4n}\cdot\frac{1}{4n}\cdot\frac{R^{4n}}{(1+R^2)^{2n}}
  =\frac{S_{4n}}{4n}.
     \end{aligned}
  \end{equation*}
  where $S_{4n}=4n\frac{\pi^{2n}}{(2n)!}$.
Thus $(\ref{LLLL})$ follows.
\end{proof}
  \begin{prop} \label{prop:real-form}
  Suppose that $\Omega \subseteq {\mathbb{H}}^n$ is a domain and $B(a,R)\Subset \Omega$ for some $R > 0$. For $u \in QSH_m(\Omega)$ and   $0< r < R$, denote
\begin{equation}\label{eq:sigma}
  \begin{aligned}
   \sigma(a,r)=\int_{B(a,r)}\Delta u \wedge \beta_n ^{n-1}.
\end{aligned}
  \end{equation}
Then, $\frac{\sigma(a,r)} {{r}^{ {\frac{4n(m-1)} m }}}$ is an increasing function of $r$ for $0< r < R$, and
\begin{equation}
v_a(u)=\lim_{r\to0} {\frac{\sigma(a,r)}{{r}^{\frac{4n(m-1)}{m}}}}
\end{equation}
exists and is
  nonnegative. It is called the $m$-Lelong number of $u$ at $a$.
\end{prop}
\begin{proof}:
For $0< r_1< r_2 < R$,  consider
  \begin{equation*}
  \begin{aligned}
v_a(r_1,r_2):=\int_{r_1<|q|\leq r_2}\Delta u \wedge \left(\Delta K_{m }\right) ^{m-1}\wedge \beta_n ^{n-m}.\\
   \end{aligned}
  \end{equation*}
Since $ K_{m }\in QSH_m$, the integrant in \eqref{eq:sigma} is a nonnegative measure on   $B(a,R) $. Without loss of generality, we may assume that $a=0$. Firstly, assume  $u\in QSH_m(B(0,R))\cap C^{\infty}(B(0,R))$.
  Then we have
  \begin{equation*}
 \begin{split}
 v_a(r_1,r_2)
  &= \int_{r_1<|q|\leq r_2}  d_0\left(   d_1 K_{m }\wedge  \Delta u
  \wedge  (\Delta K_{m }  )^{m-2}\wedge \beta_n ^{n-m}\right)\\
&=\frac { \kappa_m   }{r_2^{ 2(\kappa_m+1) }}
\int_{|q|=r_2}   \left(d_1|q|^2\wedge\triangle u \wedge  (\Delta K_{m } )^{m-2}\wedge \beta_n ^{n-m}\right)_{A}\tau( \mathbf{{n}})_{A0}~
dS\\ &\quad-\frac { \kappa_m   }{r_1^{ 2(\kappa_m+1) }}\int_{|q|=r_1}  \left( d_1|q|^2\wedge\triangle u \wedge (\Delta K_{m } )\wedge \beta_n ^{n-m}\right)_{A}\tau( \mathbf{{n}})_{A0}~
dS \\
&=(8\kappa_m)^{m-1}\left(\frac{\sigma(a,r_2)}{{r_2}^{\frac{4 n(m-1)}{m}}}-\frac{\sigma(a,r_1)}{{r_1}^{\frac{4n(m-1)}{m}}}\right)> 0,
\end{split}
\end{equation*}by using Stokes theorem,  \eqref{eq:d1-K} for $\epsilon=0$,
and
  \begin{equation*}
\begin{split}&
 \frac { \kappa_m   }{r ^{ 2(\kappa_m+1)}}
\int _ {|q|=r } \left(  d_1|q|^2\triangle u \wedge  (\Delta K_{m } )^{m-2}\wedge \beta_n ^{n-m}\right)_{A}\tau( \mathbf{{n}})_{A0}~
dS \\ = &
 \frac { 8\kappa_m   }{r ^{ 2(\kappa_m+1) }}
\int _ {|q|\leq r }    \triangle u \wedge  (\Delta K_{m } )^{m-2}\wedge \beta_n ^{n-m+1}
   = \cdots
 = \left(\frac { 8\kappa_m   }{r ^{ 2(\kappa_m+1) }}\right)^{m-1}
\int _ {|q|\leq r }    \triangle u  \wedge \beta_n ^{n -1}
\end{split}
  \end{equation*}

 Now using the convergence of $u \ast{\chi_\epsilon} \downarrow u $ and $\lim\limits_{\epsilon\to0}  \Delta (u \ast{\chi_\epsilon}) \wedge \beta_n ^{n-m}\rightarrow \Delta u \wedge \beta_n ^{n-m}$ as currents of measure type, we get the result.
\end{proof}

The proof given here also simplifies the proof of the existence of the  Lelong number  for a plurisubharmonic function in \cite{wan-wang}.

\end{document}